\documentclass[10pt]{amsart}
\usepackage{amsmath, latexsym, amsthm, amsfonts,bm,amssymb} 
\usepackage{graphicx} 
\usepackage{epsfig}
\usepackage{varioref}
\usepackage{eucal}
\usepackage{natbib} 

\bibpunct{(}{)}{;}{a}{}{,} 

\theoremstyle{plain}
\newtheorem{thm}{Theorem}

\newtheorem{lem}{Lemma}

\newtheorem{cnd}{Condition}

\theoremstyle{remark}
\newtheorem{rem}{Remark}

\newcommand{\infint}{\int_{-\infty}^{\infty}}

\def\ex{{\rm {\mathbb E\,}}}
\def\var{{\rm {\mathbb Var\,}}}

\allowdisplaybreaks

\begin{document}

\title[Deconvolution for an atomic distribution]{Deconvolution for an atomic distribution: rates of convergence}

\author{Shota Gugushvili}
\address{Department of Mathematics\\
Vrije Universiteit Amsterdam\\
De Boelelaan 1081\\
1081 HV Amsterdam\\
The Netherlands}
\email{s.gugushvili@vu.nl}
\author{Bert van Es}
\address{Korteweg-de Vries Institute for Mathematics\\
Universiteit van Amsterdam\\
P.O.\ Box 94248\\
1090 GE Amsterdam\\
The Netherlands}
\email{a.j.vanes@uva.nl}
\author{Peter Spreij}
\address{Korteweg-de Vries Institute for Mathematics\\
Universiteit van Amsterdam\\
P.O.\ Box 94248\\
1090 GE Amsterdam\\
The Netherlands}
\email{spreij@uva.nl}



\begin{abstract}
Let $X_1,\ldots, X_n$ be i.i.d.\ copies of a random variable
$X=Y+Z,$ where $ X_i=Y_i+Z_i,$ and $Y_i$ and $Z_i$ are independent
and have the same distribution as $Y$ and $Z,$ respectively.
Assume that the random variables $Y_i$'s are unobservable and that
$Y=AV,$ where $A$ and $V$ are independent, $A$ has a Bernoulli
distribution with probability of success equal to $1-p$ and $V$
has a distribution function $F$ with density $f.$ Let
the random variable $Z$ have a known distribution with density
$k.$ Based on a sample $X_1,\ldots ,X_n,$ we consider the problem
of nonparametric estimation of the density $f$ and the probability
$p.$ Our estimators of $f$ and $p$ are constructed via Fourier inversion and kernel smoothing. We derive their convergence rates over suitable functional classes. By establishing in a number of cases the lower bounds for estimation of $f$ and $p$ we show that our estimators are rate-optimal in these cases.
\end{abstract}

\date{\today}

\maketitle

\section{Introduction}
\label{deconv-intro}

Let $X_1,\ldots, X_n$ be i.i.d.\ copies of a random variable
$X=Y+Z,$ where $ X_i=Y_i+Z_i,$ and $Y_i$ and $Z_i$ are independent
and have the same distribution as $Y$ and $Z,$ respectively.
Assume that the random variables $Y_i$'s are unobservable and that
$Y=AV,$ where $A$ and $V$ are independent, $A$ has a Bernoulli
distribution with probability of success equal to $1-p$ and $V$
has a distribution function $F$ with density $f.$ Furthermore, let
the random variable $Z$ have a known distribution with density
$k.$ Based on a sample $X_1,\ldots ,X_n,$ we consider the problem
of nonparametric estimation of the density $f$ and the probability
$p.$ This problem has been recently introduced in \cite{ejs} for the case when $Z$ is normally distributed and \cite{lee} for the class of more general error distributions. It is referred to as deconvolution for an atomic
distribution, which reflects the fact that the distribution of $Y$
has an atom of size $p$ at zero and that we have to reconstruct
(`deconvolve') $p$ and $f$ from the observations from the convolution structure $X=Y+Z.$
When $p$ is known to be equal to zero, i.e.\ when $Y$ has a density, the problem reduces to the
classical and much studied deconvolution problem, see e.g.\ \cite{meister} for an introduction to the latter and many recent
references.

The above problem arises in a number of practical situations. For
instance, suppose that a measurement device is used to measure
some quantity of interest. Let it have a probability of
failure to detect this quantity equal to  $p,$ in which case it renders zero.
Repetitive measurements of the quantity of interest can be
modelled by random variables $Y_i$ defined as above. Assume that
our goal is to estimate the density $f$ and the probability of
failure $p.$ If we could use the measurements $Y_i$  directly,
then when estimating $f,$ zero measurements could be discarded and
we could use the nonzero observations to base our estimator of $f$
on. The probability $p$ could be estimated by the proportion of
zero observations. However, in practice it is often the case that
some measurement error is present. This can be modelled by random
variables $Z_i$ and assuming the additive measurement error structure, in such a case the observations are
$X_i=Y_i+Z_i.$ Now notice that due to the measurement error, the
zero $Y_i$'s cannot be distinguished from the nonzero $Y_i$'s. If
we do not want to impose parametric assumptions on $f,$ the use of
nonparametric deconvolution techniques will be unavoidable when estimating $f.$

Another example comes from evolutionary biology, see Section 4
in \cite{lee}: suppose that a virus lineage is grown in a lab for a number of days in a manner
that promotes accumulation of mutations. Plaque size can be used as
a measure of viral fitness. Assume that it is measured every day and let the mutation effect on viral
fitness be defined as a change in plaque size. If a high fitness virus is used, during any time interval in terms of mutations there are only two possibilities: either 1) no mutation, or only silent mutation occurs, or 2) a deleterious mutation occurs.
Due to the fact that a silent mutation does not affect fitness, theoretically it will not change the plaque size and hence
the mutation effect is zero for the first case. Deleterious mutations on the other hand will affect the plaque
size. Since the distribution of deleterious mutation effects is usually considered to be continuous, the distribution
of mutation effects can be expressed as a mixture of a point mass at zero, which corresponds to scenario
1), and a continuous distribution, which corresponds to scenario 2). Presence of measurement errors (which can be assumed to be additive) when measuring the plaque size leads precisely to the deconvolution problem for an atomic distribution.

Deconvolution for an atomic distribution is also closely related to
empirical Bayes estimation of a mean of a high-dimensional
normally distributed  vector, see e.g.\ \cite{zhang4} for the description of the problem and many references. In more detail, let $X_i\sim N(\theta_i,1),i=1,\ldots,n$ be i.i.d., where $N(\theta_i,1)$ denotes the normal distribution with mean $\theta_i$ and variance $1,$ and suppose that based on $X_1,\ldots,X_n$ the goal is to estimate the mean vector ${\bf \theta}=(\theta_1,\ldots,\theta_n).$ This has applications e.g.\ in denoising a noisy signal or image. It is often the case that the vector ${\bf\theta}$ is sparse in some sense in that many of $\theta_i$'s are zero or close to zero. The notion of sparsity can be naturally modelled in a Bayesian way by putting independent priors $\Pi_i(dx)=p1_{[x=0]}dx+(1-p)F(dx)$ on each component $\theta_i$ of $\theta,$ where $0\leq p<1$ and $F$ is a continuous distribution function. Notice that excess of zeros among $\theta_i$'s is matched by choosing the prior $\Pi_i$ that has a point mass at zero. In the empirical Bayes approach to estimation of $\theta$ the hyperparameters $p$ and $F$ of the priors $\Pi_i$ are estimated from the data $X_1,\ldots,X_n.$ This leads precisely to the deconvolution problem for an atomic distribution.

A related problem is estimation of the proportion of non-null effects in large-scale multiple testing framework, see e.g.\ \cite{cai}. In large-scale multiple testing one is interested in simultaneous testing of a large number of hypotheses $H_1,\ldots,H_n.$ Suppose that with every hypothesis $H_i$ there is associated a corresponding test statistic $X_i.$ A popular framework for large-scale multiple testing
is the two-group random mixture model, where one assumes that each
hypothesis $H_i$ has a certain unknown probability $\pi$ of being true (the approach is empirical Bayes in its essence) and the test statistics $X_i$ are independent and are
generated from a mixture of two densities, $X_i\sim (1-\pi)f_{\rm{null}}+\pi f_{\rm{alt}}.$ Here $\pi$ (the same for all $i$) is called the probability of null effects, $f_{\rm{null}}$ is the null density and $f_{\rm{alt}}$ is the non-null density. Often $f_{\rm{null}}$ is modelled as a density of a normal distribution $N(\mu_0,\sigma_0),$ while the density $f_{\rm{alt}}$ is modelled as a Gaussian location-scale mixture
\begin{equation*}
f_{\rm{alt}}(x)=\infint\infint \frac{1}{\sigma}\varphi\left( \frac{x-\mu}{\sigma} \right)dG(\mu,\sigma),
\end{equation*}
where $\varphi$ is the standard normal density and $G$ is the mixing distribution which is assumed to be unknown. Observe that $\pi$ in this case plays a role similar to $1-p$ in the deconvolution problem for an atomic distribution. Estimation of the probability $\pi$ and the mixing distribution $G$ based on $X_1,\ldots,X_n$ leads to a problem strongly related to the deconvolution problem for an atomic distribution.

After these motivating examples we return to the deconvolution problem for an atomic distribution and move to the construction of estimators of $p$ and $f$ (our notation is as in the first paragraph of this section). Because
of a great similarity of our problem to the classical
deconvolution problem, one natural approach to estimation of $p$ and $f$ is based on the use of Fourier
inversion and kernel smoothing, cf.\ Section 2.2.1 in
\cite{meister}. In the sequel $\phi_{\xi}$ will denote the characteristic function of a random variable $\xi.$ The Fourier transform of a function $g$ will be denoted by $\phi_g.$ Suppose that $\phi_Z(t)\neq 0$ for all
$t\in{\mathbb{R}}.$ Following \cite{ejs}, we define an estimator
$p_{ng_n}$ of $p$ as
\begin{equation}
\label{png}
p_{ng_n}=\frac{g_n}{2}\int_{-1/g_n}^{1/g_n}\frac{\phi_{emp}(t)\phi_u(g_nt)}{\phi_Z(t)}dt,
\end{equation}
where a number $g_n>0$ denotes a bandwidth, $\phi_u$ is the
Fourier transform of some fixed function (a kernel) $u$ chosen beforehand and
$\phi_{emp}(t)=n^{-1}\sum_{j=1}^n e^{itX_j}$ is the empirical
characteristic function. To make the definition of $p_{ng_n}$
meaningful, we assume that $\phi_u$ has support on $[-1,1].$ This guarantees integrability of the integrand in \eqref{png}. We
also assume that $\phi_u$ is real-valued, bounded, symmetric and integrates to two. Other conditions on
$u$ will be stated in the next section. Notice that $p_{ng_n}$ is real-valued, because for its complex conjugate we have $\overline{p_{ng_n}}=p_{ng_n}.$ The heuristics behind the
definition of $p_{ng_n}$ are the same as in \cite{ejs}: using
$\phi_X(t)=\phi_{Y}(t)\phi_Z(t)$ and $\phi_Y(t)=p+(1-p)\phi_f(t),$
we have
\begin{align*}
\lim_{g_n\rightarrow 0} \frac{g_n}{2}\int_{-1/g_n}^{1/g_n}\frac{\phi_X(t)\phi_u(g_nt)}{\phi_Z(t)}dt
& =\lim_{g_n\rightarrow 0} \frac{g_n}{2}\int_{-1/g_n}^{1/g_n} \phi_Y(t)\phi_u(g_nt)dt\\
& =\lim_{g_n\rightarrow 0} \frac{g_n}{2}\int_{-1/g_n}^{1/g_n}p\phi_u(g_nt)dt\\
&+\lim_{g_n\rightarrow 0} \frac{g_n}{2}\int_{-1/g_n}^{1/g_n}(1-p)\phi_f(t)\phi_u(g_nt)dt\\
& =p,
\end{align*}
provided $\phi_f(t)$ is integrable. The last equality follows from
the dominated convergence theorem and the fact that $\phi_u$
integrates to two. Notice that this estimator coincides with the one in
\cite{lee} when $u$ is the sinc kernel, i.e.\ $u(x)=\sin (x)/(\pi x).$ The Fourier transform of this kernel is given by $\phi_u(t)=1_{[-1,1]}(t).$
In general $p_{ng_n}$ might take on negative values, even though
for large $n$ the probability of this event will be small. In any case this is of minor importance, because we can always truncate
$p_{ng_n}$ from below at zero, i.e.\ we can define an estimator of $p$ as
$p_{ng_n}^{+}=\max(0,p_{ng_n}).$ This new estimator of $p$ has
risk (quantified by the mean square error) not larger than that of
$p_{ng_n}$:
\begin{equation*}
\ex_{p,f}[(p_{ng_n}^{+}-p)^2]\leq \ex_{p,f}[( p_{ng_n}-p)^2].
\end{equation*}
\begin{rem}
\label{notation}
In order to keep our notation compact, in the sequel instead of writing the expectation under the parameter pair $(p,f)$ as $\ex_{p,f}[\cdot],$ we will simply write $\ex[\cdot].$ \qed
\end{rem}

Next we turn to the construction of an estimator of $f.$ Let
\begin{equation}
\label{ptrunc}
\hat{p}_{ng_n}=\max(-1+\epsilon_n,\min(p_{ng_n},1-\epsilon_n)),
\end{equation}
where $0<\epsilon_n<1$ and $\epsilon_n\downarrow 0$ at a suitable
rate to be specified later on.
Notice that
$|\hat{p}_{ng_n}|\leq 1-\epsilon_n.$ Truncating $p_{ng_n}$ from below at $-1+\epsilon_n$ and not at zero will make proofs of the asymptotic results for an estimator of $f$ somewhat shorter, although truncation at zero is still a valid option. As in \cite{ejs}, we propose
the following estimator of $f,$
\begin{equation}
\label{fnhg}
f_{nh_ng_n}(x)=\frac{1}{2\pi}\infint
e^{-itx}\frac{\phi_{emp}(t)-\hat{p}_{ng_n}\phi_Z(t)}{(1-\hat{p}_{ng_n})\phi_Z(t)}\phi_w(h_nt)dt,
\end{equation}
where $w$ is a kernel function with a real-valued and symmetric Fourier transform $\phi_w$ supported
on $[-1,1]$ and $h_n>0$ is a bandwidth. Notice that $\overline{f_{nh_ng_n}(x)}=f_{nh_ng_n}(x)$ and hence $f_{nh_ng_n}(x)$ is real-valued. It is clear that
${p}_{ng_n}$ is truncated to $\hat{p}_{ng_n}$ in order to
control the factor $(1-\hat{p}_{ng_n})^{-1}$ in \eqref{fnhg}. The definition of
$f_{nh_ng_n}$ is motivated by the fact that
\begin{equation*}
f(x)=\frac{1}{2\pi}\infint e^{-itx}\frac{\phi_X(t)-p\phi_Z(t)}{(1-p)\phi_Z(t)}dt,
\end{equation*}
cf.\ equation (1.2) in \cite{ejs}. Thus $f_{nh_ng_n}$ is obtained
by replacing $\phi_X$ and $p$ by their estimators and application
of appropriate regularisation determined by the kernel $w$ and bandwidth $h.$ The estimator $f_{nh_ng_n}$ essentially coincides with the one in \cite{lee} when both $u$ and $w$ are taken to be the sinc
kernels. Again, notice that with positive probability
$f_{nh_ng_n}(x)$ might become negative for some $x\in\mathbb{R},$
a little drawback often shared by kernel-type density estimators. Some correction
method can be used to remedy this drawback, for instance one can define
$f_{nh_ng_n}^{+}(x)=\max(0,f_{nh_ng_n}(x)),$ as this does not
increase the pointwise risk of the estimator. Note that this possible negativity of $f_{nh_ng_n}$ cannot be remedied only by truncating $p_{ng_n}$ from below at zero and then using this new estimator instead of $\hat{p}_{ng_n}$ in \eqref{fnhg}. Observe also that $f_{nh_ng_n}^{+}$ can be rescaled to integrate to one and thus can be turned into a probability density. An alternative correction method to turn a possibly negative density estimator into a probability density is described in \cite{hjort}. We do not pursue
these questions any further.

In the present work we assume that the distribution of $Z$ is known. In practice this is not always the case. If the distribution of $Z$ is totally unknown, then next to the sample $X_1,\ldots,X_n$ one typically will need some additional data in order to construct consistent estimators of $f$ and $p.$ For instance, the case when additional measurements on $Z,$ say $Z_1,\ldots,Z_m,$ are available in the classical deconvolution problem with a priori known $p=0$ is dealt with in \cite{johannes}. Furthermore, one can also consider the case when the distribution of $Z$ is known up to a scale parameter. The relevant papers in the classical deconvolution context are \cite{butucea} and \cite{meist}. Although conceivable in principle, extension of our results to these cases is beyond the scope of the present work.

In the rest of the paper we concentrate on asymptotics of the
estimators $p_{ng_n}$ and $f_{nh_ng_n}.$ In particular, we derive
upper bounds on the supremum of the mean square error of the
estimator $p_{ng_n}$ and the supremum of the mean integrated
square error of the estimator $f_{nh_ng_n}$ taken over an
appropriate class of the densities $f$ and an appropriate interval
for the probability $p.$ Our results complement those in
\cite{ejs}, where the asymptotic normality of the estimators
$p_{ng_n}$ and $f_{nh_ng_n}$ is established. However, the present results
are also more general, as we consider more general error
distributions, and not necessarily the normal distribution as in
\cite{ejs}. Weak consistency of the estimators \eqref{png} and
\eqref{fnhg} based on the sinc kernel has been established under
wide conditions in \cite{lee}. Here, however, we also derive
convergence rates, much in the spirit of the classical
deconvolution problems. Notice
also that the fixed parameter asymptotics of the estimators of $p$
and $f$ were studied in \cite{lee}, in particular the rate of convergence of their estimator of $f$ (but not of $p$) was derived. On the other hand, we prefer
to study asymptotics uniformly in $p$ and $f,$ since fixed parameter
statements are difficult to interpret from the asymptotic optimality point of
view in nonparametric curve estimation, see e.g.\ \cite{low} for a
discussion. Furthermore, in case of estimation of $f$ we quantify
the risk globally in terms of the mean integrated squared error and not pointwise by the
mean squared error as done in \cite{lee}. We also derive a lower risk bound for estimation of $f,$ which shows that our estimator is rate-optimal over an appropriate functional class. Our final results are lower bounds for estimation of $p.$ These lower bounds entail rate-optimality of $p_{ng_n}$ in a large class of examples.

The structure of the paper can be outlined as follows: in Section \ref{upperbounds} we state the main results of the paper. The proofs of these results are given in Section \ref{proofs}, while the Appendix contains several technical lemmas used in Section \ref{proofs}.

\section{Results}
\label{upperbounds} The classical deconvolution problems are
usually divided into two groups, ordinary smooth deconvolution
problems and supersmooth deconvolution problems, see e.g.\
\cite{fan1} or p.\ 35 in \cite{meister}. In the former case it is assumed that the
characteristic function $\phi_Z$ of a random variable $Z$ decays
to zero algebraically at plus and minus infinity (an example of
such a $Z$ is a random variable with Laplace distribution), while in
the latter case the decay is essentially exponential (for instance,
$Z$ can be a normally distributed random variable). The rate of decay of $\phi_Z$ at infinity determines smoothness of the density of $Z$ and hence the names ordinary smooth and supersmooth. Here too we will
adopt the distinction between ordinary smooth and supersmooth deconvolution problems. The ordinary smooth deconvolution
problems for an atomic distribution will be defined by the
following condition on $\phi_Z.$
\begin{cnd}
\label{ordinarysmooth}
Let $\phi_Z(t)\neq 0$ for all $t\in{\mathbb R}$ and let
\begin{equation}
\label{orsmooth} d_0 |t|^{-\beta}\leq |\phi_Z(t)| \text{\rm{ as }} |t|\rightarrow\infty,
\end{equation}
where $d_0$ and $\beta$ are some strictly positive constants. Furthermore, let $\phi_Z$ be integrable.
\end{cnd}
\begin{rem}
\label{osmrem}
Note that the assumption of integrability of $\phi_Z$ puts certain restriction on the tail behaviour of $\phi_Z$ and therefore implicitly on $\beta$ too. In particular, in order that Condition \ref{ordinarysmooth} does not lead to an empty assumption, we must have $\beta>1.$ Notice that a lower bound on the rate of decay of $\phi_Z$ as in \eqref{orsmooth} is needed in order to derive upper risk bounds for the estimators $p_{ng_n}$ and $f_{nh_ng_n},$ cf.\ p.\ 1260 in \cite{fan1} and p.\ 35 in \cite{meister}. When deriving lower bounds for estimation of $p$ and $f,$ \eqref{orsmooth} has to be further refined by adding an explicit upper bound on the rate of decay of $\phi_Z,$ see below. \qed
\end{rem}
For the supersmooth deconvolution problems for an atomic distribution we will need the following condition on $\phi_Z.$
\begin{cnd}
\label{supersmooth}
Let $\phi_Z(t)\neq 0$ for all $t\in{\mathbb R}$ and let
\begin{equation}
\label{supsmooth} d_0 |t|^{\beta_0}e^{-|t|^{\beta}/\gamma}\leq
|\phi_Z(t)|
\text{\rm{ as }} |t|\rightarrow\infty,
\end{equation}
where $\beta_0$ is some real constant and $d_0,\beta$ and $\gamma$ are some strictly positive constants. Furthermore, let $\phi_Z$ be integrable.
\end{cnd}
Next we need to impose conditions on the class of target
densities $f.$
\begin{cnd}
\label{fclass} Define the class of target densities $f$ as
\begin{equation}
\label{sobolev} {\Sigma}(\alpha,K_{\Sigma})=\left\{ f: \infint
|\phi_f(t)|^2(1+|t|^{2\alpha})dt\leq K_{\Sigma} \right\},
\end{equation}
Here $\alpha$ and $K_{\Sigma}$ are some strictly
positive numbers.
\end{cnd}
Smoothness conditions of this type are typical in nonparametric curve
estimation problems, cf.\ p.\ 25 in \cite{tsyb} or p.\ 34 in
\cite{meister}. Some smoothness assumptions have to be imposed on the class of target densities, because e.g.\ the class of all continuous densities is usually too large to be handled when dealing with uniform asymptotics.
A possibility, different from Condition \ref{fclass}, is to assume that $f$ belongs to the class of supersmooth densities
\begin{equation*}
\Sigma(\alpha,\gamma,K_{\Sigma}) = \left\{ f: \infint |\phi_f(t)|^2\exp(2\gamma|t|^{\alpha})dt\leq K_{\Sigma} \right\},
\end{equation*}
for some strictly positive $\alpha,\gamma$ and $K_{\Sigma}.$ The class $\Sigma(\alpha,\gamma,K_{\Sigma})$ is much smaller than the class ${\Sigma}(\alpha,K_{\Sigma})$ and the estimators $p_{ng_n}$ and $f_{ng_nh_n}$ will enjoy better convergence rates in this case than in the case when the class of target densities is ${\Sigma}(\alpha,K_{\Sigma}),$ cf.\ \cite{tsyb2} and \cite{tsyb1} for a similar result in the classical deconvolution problem. In order not to overstretch the length of the paper, we decided however not to cover this case in the present work.

\begin{rem}
\label{lesssim}
In the sequel we will use the symbols $\lesssim$ and $\gtrsim$ to compare two sequences $a_n$ and $b_n$ indexed by $n,$ meaning respectively that $a_n$ is less or equal than $b_n$ for all $n,$ or greater or equal, up to a universal constant that does not depend on $n.$ \qed
\end{rem}

The following theorem deals with asymptotics of the estimator
$p_{ng_n}.$ Its proof, as well as the proofs of all other results
of the paper, is given in Section \ref{proofs}.
\begin{thm}
\label{thmpng} Let a function $u$ be such that its Fourier transform $\phi_u$ is symmetric, real-valued,
continuous in some neighbourhood of zero and is supported on $[-1,1].$ Furthermore, let
\begin{equation}
\label{u} \int_{-1}^1 \phi_u(t)dt=2, \quad \left|
\frac{\phi_u(t)}{t^{\alpha}} \right|\leq U{\text{ for all }}t\in\mathbb{R},
\end{equation}
where the constant $\alpha$ is the same as in Condition
\ref{fclass}, $U$ is a strictly positive constant and for $t=0$ the ratio $\phi_u(t)t^{-\alpha}$ is defined by
continuity at zero as the limit $\lim_{t\rightarrow 0}\phi_u(t)t^{-\alpha},$ which we assume to exist. Then

(i) under Condition \ref{ordinarysmooth}, by selecting $g_n = d
n^{-1/(2\alpha+2\beta)}$ for some constant $d>0,$ we have
\begin{equation}
\label{mseord2}
\sup_{f\in\Sigma(\alpha,K_{\Sigma}),p\in[0,1)}\ex
[(p_{ng_n}-p)^2]\lesssim n^{-(2\alpha+1)/(2\alpha+2\beta)};
\end{equation}

(ii) under Condition \ref{supersmooth}, by selecting $g_n
=(4/\gamma)^{1/\beta} (\log n)^{-1/\beta},$ we have
\begin{equation}
\label{msesuper2}
\sup_{f\in\Sigma(\alpha,K_{\Sigma}),p\in[0,1)}\ex
[(p_{ng_n}-p)^2]\lesssim (\log n)^{-(2\alpha+1)/\beta}.
\end{equation}
\end{thm}

Thus the rate of convergence of the estimator $p_{ng_n}$ is slower
than the root-${n}$ rate for estimation of a finite-dimensional
parameter in regular parametric models. For Theorem \ref{thmpng} (ii) this is evident, while for Theorem \ref{thmpng} (i) this follows from Remark \ref{osmrem}, which entails the fact that $2\alpha+1<2\alpha+2\beta.$ However, see Theorems \ref{lbnd-p} and \ref{lbnd-p-os} below, where for a practically important case of a normally distributed $Z,$ as well as $Z$ with ordinary smooth distribution, by establishing the lower bounds for estimation of $p$ we show that the slow convergence rate is intrinsic to the deconvolution problem and is not a quirk of our particular estimator.

\begin{rem}
\label{kerrem}
The function $u$ in the statement of Theorem \ref{thmpng} will not be a probability density, not even a function that integrates to one, and hence by calling it a kernel we somewhat abuse the established terminology in kernel estimation. Notice that condition \eqref{u} and the assumption $\alpha>0$ in Condition \ref{fclass} preclude the kernel $u$ from being the sinc kernel. We refer to \cite{ejs} for one particular example of $u$ that produced good results in simulations. Its Fourier transform is given by
\begin{equation*}
\phi_u(t)=\frac{693}{8}t^6(1-t^2)^2 1_{[-1,1]}(t).
\end{equation*}
Here $\alpha=6$ and $U=693/8.$ An explicit, but rather complicated expression for $u$ can be found in \cite{ejs}. \qed
\end{rem}

Next we will study the asymptotic behaviour of the estimator
$f_{nh_ng_n}$ of $f.$ We select the mean integrated square error
as a criterion of its performance.

Due to technical reasons, see the proof of Theorem \ref{thmfnhg}, in the ordinary smooth case it is convenient to split the sample $X_1,\ldots,X_n$ into two parts and next to base the estimator $\hat{p}_{ng_n}$ on the first part of the sample only, i.e.\ on $X_1,\ldots,X_{\lfloor n/2\rfloor},$ and to redefine $f_{nh_ng_n}$ as
\begin{equation}
\label{ssplf}
f_{nh_ng_n}(x)=\frac{1}{2\pi} \infint e^{-itx} \frac{ \widetilde{\phi}_{emp}(t) - \hat{p}_{ng_n} }{(1-\hat{p}_{ng_n}) \phi_Z(t)} \phi_w(h_n t)dt,
\end{equation}
where
\begin{equation*}
\widetilde{\phi}_{emp}(t) = \frac{1}{n-\lfloor n/2 \rfloor} \sum_{j= \lfloor n/2 \rfloor +1 }^n e^{itX_j}.
\end{equation*}
Thus $\widetilde{\phi}_{emp}$ is based on the second half of the sample $X_1,\ldots,X_n$ only. Note that $\ex[\phi_{emp}(t)]=\ex[\widetilde{\phi}_{emp}(t)]= \phi_X(t).$ From now on we will assume that $p_{ng_n}$ and $f_{nh_ng_n}$ are defined in this way in the ordinary smooth case, but will retain the old definition in the supersmooth case. Splitting the sample does not affect the convergence rate of $f_{nh_ng_n}$ in the ordinary smooth case, but only the constant factor in the upper bound on its mean integrated squared error. The general case without sample splitting in principle can also be handled, but we anticipate longer and more technical proofs, cf. the remarks at the end of the proof of Theorem \ref{thmfnhg}. Since in the present work we are only concerned with convergence rates, sample splitting does not lead to a significant loss of generality.

The following theorem holds.
\begin{thm}
\label{thmfnhg} Let a kernel $u$ satisfy the
assumptions in Theorem \ref{thmpng}. Furthermore, let a kernel $w$
be such that its Fourier transform is symmetric, real-valued, is supported on $[-1,1]$ and
\begin{equation}
\label{k}
\phi_w(0)=1, \quad |\phi_w(t)-1|\leq W|t|^{\alpha}\text{ for all }t\in\mathbb{R},
\quad \int_{-1}^1 |\phi_w(t)|^2dt<\infty,
\end{equation}
where $W$ is some strictly positive constant. Moreover, let $p\in[0,p^*],$ where
$p^*<1.$ Then

(i) under Condition \ref{ordinarysmooth}, by selecting $h_n
=d (n-\lfloor n/2\rfloor)^{-1/(2\alpha+2\beta+1)}$ for some $d>0,$ $g_n=d \lfloor n/2 \rfloor^{-1/(2\alpha+2\beta)}$ and $\epsilon_n=(\log 3n)^{-1},$ we have
\begin{equation}
\label{miseord2}
\sup_{f\in\Sigma(\alpha,K_{\Sigma}),p\in[0,p^*]}\ex \left[\infint(f_{nh_ng_n}(x)-f(x))^2dx\right]\lesssim n^{-2\alpha/(2\alpha+2\beta+1)},
\end{equation}
where $f_{nh_ng_n}$ is defined by \eqref{ssplf}.

(ii) under Condition \ref{supersmooth}, by selecting $h_n = g_n =
(4/\gamma)^{1/\beta} (\log n)^{-1/\beta}$ and $\epsilon_n=(\log 3n)^{-1},$ we have
\begin{equation}
\label{misesuper2}
\sup_{f\in\Sigma(\alpha,K_{\Sigma}),p\in[0,p^*]}\ex \left[\infint(f_{nh_ng_n}(x)-f(x))^2dx\right]\lesssim  (\log n)^{-2\alpha/\beta},
\end{equation}
where $f_{nh_ng_n}$ is defined by \eqref{fnhg}.
\end{thm}

\begin{rem}
\label{pst}
As it will become clear from the proof of this theorem, without the assumption $p^*<1$ one cannot study the asymptotics of $f_{nh_ng_n}$ uniformly in $(p,f)$ for $p\in[0,p^*]$ and $f\in\Sigma(\alpha,K_{\Sigma}).$ Since $p^*$ is allowed to be arbitrarily close to $1,$ from a practical point of view $p^*<1$ is not an important restriction. Observe that one can also study the case when $p^*=p^*_n$ depends on the sample size $n$ and $p^*_n\rightarrow 1$ at a suitable rate.\qed
\end{rem}

\begin{rem}
\label{pbch}
The condition $h_n=g_n$ in Theorem \ref{thmfnhg} (ii) is imposed for simplicity of the proofs
only. In practice the two bandwidths need not be the same, cf.\
\cite{ejs}, where unequal $h_n$ and $g_n$ are used in simulation
examples. Also notice that our conditions on $h_n$ and $g_n$ in Theorems \ref{thmpng} and \ref{thmfnhg} are
of asymptotic nature. For practical suggestions
on bandwidth selection for the case when both $u$ and $w$ are
sinc kernels, see \cite{lee}, where also a number of simulation
examples is given. \qed
\end{rem}

\begin{rem}
\label{kerwrem}
We refer to \cite{ejs} for one particular example of a kernel $w.$ Any kernel that is known to produce good results in the classical deconvolution problem can be used as a kernel $w.$ A relevant paper on the choice of a kernel in the context of the classical deconvolution problems is \cite{delaiglehall}, to which we refer for a discussion and more examples. \qed
\end{rem}

The upper risk bounds derived in Theorem \ref{thmfnhg} coincide with the upper risk bounds for kernel-type estimators in the classical deconvolution problems, i.e.\ in the case when $p$ is a priori known to be zero, see Theorem 2.9 in \citet{meister}. Naturally, a discussion on the optimality of convergence rates
of the estimators $f_{nh_ng_n}$ and $p_{ng_n}$ is in order. Let
$\widetilde{f}_n$ denote an arbitrary estimator of $f$ based on a
sample $X_1,\ldots,X_n.$ Consider
\begin{equation*}
\mathcal{R}_n^{*}\equiv \inf_{\widetilde{f}_n}\sup_{f\in\Sigma,p\in[0,p^*]} \ex\left[\infint(\tilde{f}_{n}(x)-f(x))^2dx\right],
\end{equation*}
i.e.\ the minimax risk for estimation of $f$ over some functional class $\Sigma$ and the interval $[0,p^{\ast}]$ for $p$ that is associated with our statistical model, cf.\ p.\ 78 in \cite{tsyb}. Notice that
\begin{equation*}
\mathcal{R}_n^{*}\geq \inf_{\widetilde{f}_n}\sup_{f\in\Sigma,p=0} \ex\left[\infint(\tilde{f}_{n}(x)-f(x))^2dx\right].
\end{equation*}
The quantity on the right-hand side coincides with the minimax
risk for estimation of a density $f$ in the classical
deconvolution problem, i.e.\ when $p=0$ and the random variable $Y$ has a
density $f$. Using this fact, by Theorem 2.14 of \cite{meister}
it is easy to obtain lower bounds for $\mathcal{R}_n^{*},$ but first we need to formulate two addition conditions on the rate of decay of $\phi_Z$ at plus and minus infinity. These two conditions correspond to the ordinary smooth and supersmooth deconvolution problems, cf.\ Conditions \ref{ordinarysmooth} and \ref{supersmooth}.

\begin{cnd}
\label{osmlb}
Let $\phi_Z$ be such that
\begin{equation*}
|\phi_Z(t)|\leq \frac{d_1}{1+|t|^{\beta}}, \quad |\phi_Z^{\prime}(t)|\leq \frac{d_1}{1+|t|^{\beta}} \quad \text{\rm{for all }} t\in\mathbb{R}
\end{equation*}
for some strictly positive constants $d_1$ and $\beta.$
\end{cnd}

\begin{cnd}
\label{spsmlb}
Let $\phi_Z$ be such that
\begin{equation*}
|\phi_Z(t)|\leq d_1 e^{-|t|^{\beta}/\gamma}, \quad |\phi_Z^{\prime}(t)|\leq d_1 e^{-|t|^{\beta}/\gamma} \quad \text{\rm{for all }} t\in\mathbb{R}
\end{equation*}
for some strictly positive constants $d_1,\beta$ and $\gamma.$
\end{cnd}

The following result holds.

\begin{thm}
\label{lbnd-f} Let $\widetilde{f}_n$ denote any estimator of $f$ based on a sample $X_1,\ldots,X_n$ and let $\alpha\geq 1/2.$ Suppose that $K_{\Sigma}$ is large enough. Then

(i) under Condition \ref{osmlb} we have
\begin{equation}
\label{lbndf1}
\inf_{\widetilde{f}_n}\sup_{f\in\Sigma(\alpha,K_{\Sigma}),p\in[0,p^*]}\ex \left[\infint(\hat{f}(x)-f(x))^2dx\right]\gtrsim n^{-2\alpha/(2\alpha+2\beta+1)};
\end{equation}

(ii) under Condition \ref{spsmlb} the inequality
\begin{equation}
\label{lbndf2}
\inf_{\widetilde{f}_n}\sup_{f\in\Sigma(\alpha,K_{\Sigma}),p\in[0,p^*]}\ex \left[\infint(\hat{f}(x)-f(x))^2dx\right]\gtrsim
(\log n)^{-2\alpha/\beta}
\end{equation}
holds.
\end{thm}

These lower bounds are of the
same order as upper bounds in Theorem \ref{thmfnhg}. It then follows that our estimator of $f$ is rate-optimal under the combined conditions in Theorems \ref{thmfnhg} and \ref{lbnd-f}. For a discussion on the conditions in Theorem 3 see p.\ 35 in \citet{meister}.

Derivation of the lower risk bounds for estimation of probability $p$ appears to be more involved. We will establish the lower bound for the case when $Z$ follows the standard normal distribution. This is an important case, as the assumption of normality of measurement errors is frequently imposed in practice. The following result holds true.

\begin{thm}
\label{lbnd-p} Let $Z$ have the standard normal distribution and let $\widetilde{p}_n$ denote any estimator of $p$ based on a sample $X_1,\ldots,X_n.$ Then
\begin{equation}
\label{lbndp}
\inf_{\widetilde{p}_n}\sup_{f\in\Sigma(\alpha,K_{\Sigma}),p\in[0,1)}\ex \left[(\widetilde{p}_n-p)^2\right]\gtrsim
(\log n)^{-(\alpha+1/2)}
\end{equation}
holds.
\end{thm}

A consequence of this theorem and \eqref{msesuper2} is that our estimator $p_{ng_n}$ is rate-optimal in the case when $Z$ follows the normal distribution. 

The arguments used in the proof of Theorem \ref{lbnd-p} can be easily extended to the case when the distribution of $Z$ is ordinary smooth. Below we provide the corresponding statement in the ordinary smooth case.
\begin{thm}
\label{lbnd-p-os} Let the characteristic function of $Z$ satisfy Condition \ref{osmlb}
for $\beta>1/2.$
Let $\widetilde{p}_n$ denote any estimator of $p$ based on the sample $X_1,\ldots,X_n.$ Then
\begin{equation*}
\inf_{\widetilde{p}_n}\sup_{f\in\Sigma(\alpha,K_{\Sigma}),p\in[0,1)}\ex \left[(\widetilde{p}_n-p)^2\right]\gtrsim
 n^{-(2\alpha+1)/(2\alpha+2\beta)}
\end{equation*}
holds.
\end{thm}
This theorem and Theorem \ref{png} (i) imply that under the combined conditions in Theorems \ref{thmpng} (i) and \ref{lbnd-p-os} the estimator $p_{ng_n}$ is rate-optimal.

\section{Proofs}
\label{proofs}
\begin{proof}[Proof of Theorem \ref{thmpng}] The proof uses some arguments from \cite{fan1}. To make the notation less cumbersome, let
$\sup_{f,p}\equiv\sup_{f\in\Sigma(\alpha,K_{\Sigma}),p\in[0,1)}.$
We first prove (i). We have
\begin{equation}
\label{biasvar}
\sup_{f,p}
\ex[(p_{ng_n}-p)^2] \leq
\sup_{f,p}(\ex[p_{ng_n}]-p)^2+\sup_{f,p}\var[p_{ng_n}].
\end{equation}
Observe that
\begin{equation}
\label{starp7}
\begin{split}
|\ex[p_{ng_n}]-p|&=\frac{1-p}{2}\left|\int_{-1}^1\phi_f\left(
\frac{t}{g_n} \right)\phi_u(t)dt\right|\\
& \leq \frac{1}{2} \int_{-1}^1 \left| \phi_f\left( \frac{t}{g_n} \right) \left(\frac{t}{g_n}\right)^{\alpha} \right| \left| \left(\frac{g_n}{t}\right)^{\alpha} \phi_u(t) \right| 1_{[t\neq 0]} dt\\
&\leq \frac{1}{2} \sqrt{ \int_{-1}^1 \left| \phi_f\left( \frac{t}{g_n} \right) \left(\frac{t}{g_n}\right)^{\alpha} \right|^2 dt} \sqrt{ \int_{-1}^1 \left| \frac{\phi_u(t)}{t^{\alpha}} g_n^{\alpha} \right|^2 1_{[t\neq 0]} dt }\\
&\leq \frac{1}{\sqrt{2}}
\sqrt{K_{\Sigma}}U g_n^{\alpha+1/2},
\end{split}
\end{equation}
where we used \eqref{u}, \eqref{sobolev} and the Cauchy-Schwarz inequality. Therefore
\begin{equation}
\label{bias1} \sup_{f,p}(\ex[p_{ng_n}]-p)^2\lesssim
g_n^{2\alpha+1}
\end{equation}
holds. Furthermore, using independence of the random variables $X_i$'s,
\begin{align*}
\var[p_{ng_n}] & =\frac{1}{4}\frac{g^2_n}{n}\var \left[ \int_{-1/g_n}^{1/g_n} e^{itX_1} \frac{\phi_u(g_n t)}{\phi_Z(t)}dt \right]\\
& \leq \frac{1}{4}\frac{g^2_n}{n} \ex\left[ \left( \int_{-1/g_n}^{1/g_n} e^{itX_1} \frac{\phi_u(g_n t)}{\phi_Z(t)}dt \right)^2 \right]\\
& = \frac{1}{4}\frac{g^2_n}{n} \int_{-\infty}^{\infty} \left( \int_{-1/g_n}^{1/g_n} e^{itx} \frac{\phi_u(g_n t)}{\phi_Z(t)}dt \right)^2 q(x)dx,
\end{align*}
where $q$ is the density of $X_1.$ Notice that
\begin{equation*}
q(x)= \frac{1}{2\pi} \int_{-\infty}^{\infty} e^{-itx} \phi_Y(t)\phi_Z(t)dt \leq \frac{1}{2\pi}  \int_{-\infty}^{\infty} |\phi_Z(t)|dt<\infty,
\end{equation*}
where we used integrability of $\phi_Z.$
Therefore
\begin{align*}
\var[p_{ng_n}] & \lesssim \frac{g_n^2}{n} \int_{-\infty}^{\infty} \left( \int_{-1/g_n}^{1/g_n} e^{itx} \frac{\phi_u(g_n t)}{\phi_Z(t)}dt \right)^2 dx\\
& = \frac{1}{2\pi} \frac{g_n^2}{n} \int_{-\infty}^{\infty} \left| \frac{\phi_u(g_n t)}{\phi_Z(-t)} \right|^2 dt\\
& = \frac{1}{2\pi} \frac{g_n}{n} \int_{-1}^{1} \left| \frac{\phi_u(t)}{\phi_Z(-t/g_n)} \right|^2 dt
\end{align*}
by Parseval's identity. This inequality and an argument as on p.\ 1266 of \cite{fan1} entail
that
\begin{equation}
\label{var5}
\sup_{f,p}\var[p_{ng_n}]\lesssim
\frac{1}{ng_n^{2\beta-1}}.
\end{equation}
Formula \eqref{mseord2} is then a consequence of \eqref{biasvar},
\eqref{bias1}, \eqref{var5} and our specific choice of $g_n$ in
(i).

Now we prove (ii). Since the first term on the right-hand side of
\eqref{biasvar} can be treated as in the ordinary smooth case (in
particular \eqref{bias1} holds), we concentrate on the second
term. Using independence of the random variables $X_i$'s,
\begin{align}
\label{var1}
\var[p_{ng_n}] & =\frac{1}{4}\frac{1}{n}\var \left[ \int_{-1}^1 e^{itX_1/g_n} \frac{\phi_u(t)}{\phi_Z(t/g_n)}dt \right]\nonumber\\
& \leq \frac{1}{4}\frac{1}{n} \left( \int_{-1}^1 \left| \frac{\phi_u(t)}{\phi_Z(t/g_n)} \right| dt \right)^2.
\end{align}
By the same arguments as on pp.\ 1265--1266 of \cite{fan1}, one can show
that
\begin{equation}
\label{2star} \int_{-1}^1 \left| \frac{\phi_u(t)}{\phi_Z(t/g_n)}
\right| dt\leq
\begin{cases}
C^{'} e^{1/(\gamma g_n^{\beta})},  & \mbox{if }\beta_0\geq 0\\
C^{'} g_n^{\beta_0} e^{1/(\gamma g_n^{\beta})}, & \mbox{if
}\beta_0< 0,
\end{cases}
\end{equation}
where the constant $C^{'}$ does not depend on $n.$ In either case,
because of our choice of $g_n,$ the righthand side of
\eqref{2star} is of order $o(n^{1/3}).$ This and \eqref{var1} imply that
\begin{equation*}
\sup_{f,p}\var[p_{ng_n}]=o(n^{-1/3}).
\end{equation*}
The latter together with \eqref{biasvar}, \eqref{bias1} and our choice of $g_n$ in (ii) proves
\eqref{msesuper2}.
\end{proof}
\begin{proof}[Proof of Theorem \ref{thmfnhg}]
We use the shorthand notation
$\sup_{f,p}\equiv\sup_{f\in\Sigma(\alpha,K_{\Sigma}),p\in[0,p^*]}.$
By Fubini's theorem and the standard squared bias plus variance decomposition we have
\begin{align*}
\sup_{f,p}\ex\left[\infint(f_{nh_ng_n}(x)-f(x))^2dx\right]&\leq \sup_{f,p} \infint(\ex[f_{nh_ng_n}(x)]-f(x))^2dx\\
&+\sup_{f,p}\infint\var[f_{nh_ng_n}(x)]dx\\
&=T_1+T_2.
\end{align*}
Keeping in mind the remarks surrounding \eqref{ssplf}, let
\begin{equation*}
\hat{f}_{nh_n}(x)=\frac{1}{2\pi}\infint e^{-itx}\frac{\widetilde{\phi}_{emp}(t)\phi_w(h_nt)}{\phi_Z(t)}dt
\end{equation*}
in the ordinary smooth case, while
\begin{equation*}
\hat{f}_{nh_n}(x)=\frac{1}{2\pi}\infint e^{-itx}\frac{{\phi}_{emp}(t)\phi_w(h_nt)}{\phi_Z(t)}dt
\end{equation*}
in the supersmooth case. Introduce
\begin{equation}
\label{nomeri**}
f_{nh_n}(x)=\frac{\hat{f}_{nh_n}(x)}{1-p}-\frac{p}{1-p}w_{h_n}(x),
\end{equation}
where $w_{h_n}(x)=(1/h_n)w(x/h_n).$ We first study $T_1,$ i.e.\
the supremum of the integrated squared bias. By the
$c_2$-inequality it can be bounded as
\begin{align*}
T_1&\lesssim\sup_{f,p}\infint(\ex[f_{nh_n}(x)]-f(x))^2dx\\
&+\sup_{f,p}\infint(\ex[f_{nh_ng_n}(x)-f_{nh_n}(x)])^2dx\\
&=T_3+T_4.
\end{align*}
By Parseval's identity and the dominated convergence theorem
\begin{align*}
\infint(\ex[f_{nh_n}(x)]-f(x))^2dx&=\frac{1}{2\pi}\infint |\phi_f(t)|^2|\phi_w(h_nt)-1|^2dt\\
&=h_n^{2\alpha}\frac{1}{2\pi}\infint |t|^{2\alpha}|\phi_f(t)|^2\frac{|\phi_w(h_nt)-1|^2}{|h_nt|^{2\alpha}} 1_{[t\neq 0]} dt\\
&\lesssim h^{2\alpha}_n.
\end{align*}
Here in the second equality we used the fact that $\phi_w(0)=1.$ The dominated convergence theorem is applicable because of
Condition \ref{fclass} and \eqref{k}. Hence $T_3\lesssim
h_n^{2\alpha}$ in view of the fact that
$f\in\Sigma(\alpha,K_{\Sigma}).$ It is also straightforward to see that in fact $\sup_{f,p}T_3\lesssim
h_n^{2\alpha}.$ We deal with
$T_4.$ By the $c_2$-inequality
\begin{align*}
\infint(\ex[f_{nh_ng_n}(x)-f_{nh_n}(x)])^2dx&\lesssim \left(\ex\left[ \frac{\hat{p}_{ng_n}-p}{(1-\hat{p}_{ng_n})(1-p)}\right] \right)^2\infint (w_{h_n}(x))^2dx\\
&+\infint \left(\ex\left[ \hat{f}_{nh_n}(x)\frac{(\hat{p}_{ng_n}-p)}{(1-\hat{p}_{ng_n})(1-p)} \right]\right)^2dx\\
&=T_5+T_6.
\end{align*}
Notice that
\begin{equation*}
\infint (w_{h_n}(x))^2dx=\frac{1}{h_n} \infint (w(x))^2dx <\infty,
\end{equation*}
because by our assumptions and Parseval's identity $w$ is square integrable. We first consider $T_5.$ By the Cauchy-Schwarz inequality we have
\begin{equation*}
T_5 \leq \frac{1}{h_n}\infint (w(u))^2du \ex\left[ \frac{(\hat{p}_{ng_n}-p)^2}{(1-\hat{p}_{ng_n})^2(1-p)^2}\right].
\end{equation*}
With our choice of the smoothing parameters $h_n$ and $g_n$ it
follows from Lemma \ref{lemma2} of the Appendix that
$\sup_{p,f}T_5\lesssim g_n^{2\alpha}.$ Now let us
turn to $T_6.$ By the Cauchy-Schwarz inequality
\begin{equation*}
T_6\leq
\ex\left[\frac{(\hat{p}_{ng_n}-p)^2}{(1-\hat{p}_{ng_n})^2(1-p)^2}\right]
\infint \ex[ (\hat{f}_{nh_n}(x))^2]dx.
\end{equation*}
By Lemma \ref{lemma2} of the Appendix the first term in the
product in the above display is of order
$g_n^{2\alpha+1}.$ The same holds true for its
supremum over $f$ and $p.$ Hence it remains to study the second
factor in the above upper bound on $T_6.$ We have
\begin{align*}
\infint \ex[ (\hat{f}_{nh_n}(x))^2]dx&=\infint \var[
\hat{f}_{nh_n}(x)]dx+\infint
(\ex[\hat{f}_{nh_n}(x)])^2dx\\
&=T_7+T_8.
\end{align*}
Let the function $W_n$ is defined by
\begin{equation*}
W_n(x)=\frac{1}{2\pi}\int_{-1}^1
e^{-itx}\frac{\phi_w(t)}{\phi_Z(t/h_n)}dt.
\end{equation*}
Notice that by independence of $X_i$'s
\begin{equation*}
T_7=\frac{1}{nh_n^2}\infint \var\left[
W_n\left(\frac{x-X_1}{h_n}\right) \right]dx\leq
\frac{1}{nh_n^2}\infint \ex\left[\left(
W_n\left(\frac{x-X_1}{h_n}\right)\right)^2 \right]dx
\end{equation*}
in the supersmooth case, and
\begin{equation*}
T_7\leq
\frac{1}{(n-\lfloor n/2\rfloor)h_n^2}\infint \ex\left[\left(
W_n\left(\frac{x-X_1}{h_n}\right)\right)^2 \right]dx
\end{equation*}
in the ordinary smooth case.
Then by Fubini's theorem
\begin{align*}
T_7&\leq\frac{1}{nh_n^2} \infint \infint \left(
W_n\left(\frac{x-s}{h_n}\right) \right)^2q(s)dsdx\\
&=\frac{1}{nh_n^2} \infint \infint \left(
W_n\left(\frac{x-s}{h_n}\right) \right)^2dxq(s)ds\\
&=\frac{1}{nh_n}\infint\infint (W_n(x))^2dx q(s)ds\\
&=\frac{1}{nh_n} \int_{-1}^1
\frac{|\phi_w(t)|^2}{|\phi_Z(t/h_n)|^2}dt
\end{align*}
in the supersmooth case, and
\begin{equation*}
T_7\leq\frac{1}{(n-\lfloor n/2 \rfloor)h_n} \int_{-1}^1
\frac{|\phi_w(t)|^2}{|\phi_Z(t/h_n)|^2}dt
\end{equation*}
in the ordinary smooth case. Here we used the fact that $q,$ being a probability density,
integrates to one, as well as Parseval's identity. The integrals in
the last equalities of the above two displayed formulae can be analysed by exactly
the same arguments as on pp.\ 1265-1266 in \citet{fan1}. Thus
\begin{equation}
\label{vii} T_7 \lesssim
\begin{cases}
\frac{1}{nh_n^{2\beta+1}}, & \mbox{if }Z \mbox{ is ordinary
smooth},\\
\frac{1}{nh_n} e^{2/(\gamma h_n^{\beta})}, & \mbox{if
}Z \mbox{ is supersmooth and }\beta_0\geq 0,\\
\frac{h_n^{2\beta_0-1}}{n} e^{2/(\gamma h_n^{\beta})}, & \mbox{if
}Z \mbox{ is supersmooth and }\beta_0<0.
\end{cases}
\end{equation}
The same order bounds hold for
$\sup_{f,p}T_7$ as well. As a consequence, $\sup_{f,p}T_7\rightarrow 0.$ Let us now study $T_8.$ By Parseval's identity
and the fact that $|\phi_Y(t)|\leq 1,$ we have
\begin{align*}
T_8&=\infint \left( \frac{1}{2\pi} \int_{-1/h_n}^{1/h_n} e^{-itx}\phi_Y(t)\phi_w(h_nt)dt \right)^2dx\\
&=\frac{1}{2\pi}\infint |\phi_Y(t)\phi_w(h_nt)|^2 1_{[-h^{-1},h^{-1}]}(t) dt\\
&\leq \frac{1}{h_n}\frac{1}{2\pi}\int_{-1}^1 |\phi_w(t)|^2dt\\
&\lesssim\frac{1}{h_n},
\end{align*}
where the last line follows from our assumptions on $w.$
It follows that $\sup_{p,f} T_8 \lesssim 1/h_n.$ Combination of the above bounds on $\sup_{p,f}T_7$ and $\sup_{p,f}T_8$
entails that $\sup_{f,p}T_6\lesssim g_{n}^{2\alpha},$ where we also used the fact that $g_n\lesssim h_n.$ Therefore $T_4,$ as well as
$T_1,$ i.e.\ the supremum of the integrated squared bias, is of
order $h_{n}^{2\alpha}.$ For the ordinary smooth case this gives
an upper bound of order $n^{-2\alpha/(2\alpha+2\beta+1)}$ on $T_1,$ while
for the supersmooth case an upper bound of order $(\log
n)^{-2\alpha/\beta}.$

Now we turn to $T_2,$ i.e.\ the supremum of the integrated
variance. We have
\begin{align*}
\infint \var[f_{nh_ng_n}(x)]dx&=\infint\var[f_{nh_ng_n}(x)-f_{nh_n}(x)+f_{nh_n}(x)]dx\\
&\lesssim\infint\var[f_{nh_n}(x)]dx+\infint\var[f_{nh_ng_n}(x)-f_{nh_n}(x)]dx\\
&=T_9+T_{10},
\end{align*}
where we used the fact that for random variables $\xi$ and $\eta$
\begin{equation*}
\var[\xi+\eta] \leq 2(\var[\xi]+\var[\eta]).
\end{equation*}
Since $T_9$ up to a constant is the same as $T_7,$ cf.\ \eqref{nomeri**}, the term $\sup_{f,p}T_9$ can be bounded
as before, see \eqref{vii}. We consider $T_{10}.$ Let
$\psi_n$ be as in \eqref{psin} in the proof of Lemma \ref{lemma2} of the Appendix. Then
\begin{align*}
T_{10}&\leq  \infint\ex[(f_{nh_ng_n}(x)-f_{nh_n}(x))^2 1_{[|\hat{p}_{ng_n}-p|>\psi_n]}]dx\\
&+\infint\ex[(f_{nh_ng_n}(x)-f_{nh_n}(x))^2 1_{[|\hat{p}_{ng_n}-p|\leq\psi_n]}]dx\\
&=T_{11}+T_{12}.
\end{align*}
By the $c_2$-inequality
\begin{align*}
T_{11} & \lesssim \frac{1}{h_n}\infint (w(x))^2dx\ex\left[\frac{(\hat{p}_{ng_n}-p)^2}{(1-\hat{p}_{ng_n})^2(1-p)^2} 1_{[|\hat{p}_{ng_n}-p|>\psi_n]} \right]\\
&+\infint\ex\left[(\hat{f}_{nh_n}(x))^2\frac{(\hat{p}_{ng_n}-p)^2}{(1-\hat{p}_{ng_n})^2(1-p)^2} 1_{[|\hat{p}_{ng_n}-p|>\psi_n]}\right]dx\\
&=T_{13}+T_{14}.
\end{align*}
Since $T_{13}\lesssim
h_n^{-1}\epsilon_n^{-2}\sup_{p,f}P(|\hat{p}_{ng_n}-p|>\psi_n),$
which follows from the fact that
\begin{equation*}
\frac{(\hat{p}_{ng_n}-p)^2}{(1-\hat{p}_{ng_n})^2(1-p)^2} \leq \frac{2(1-\epsilon_n)^2+2p^{*2}}{\epsilon_n^2(1-p^*)^2},
\end{equation*}
by Lemma \ref{lemma3} of the Appendix with our conditions on $h_n$ and $\epsilon_n$ it certainly holds true that $\sup_{p,f}T_{13}\lesssim h_{n}^{2\alpha}.$ As far
as $T_{14}$ is concerned, by Fubini's theorem and Parseval's
identity
\begin{align*}
T_{14}&=\ex\left[ \frac{(\hat{p}_{ng_n}-p)^2}{(1-\hat{p}_{ng_n})^2(1-p)^2} 1_{[|\hat{p}_{ng_n}-p|>\psi_n]} \infint (\hat{f}_{nh_n}(x))^2dx \right]\\
&=\ex\left[ \frac{(\hat{p}_{ng_n}-p)^2}{(1-\hat{p}_{ng_n})^2(1-p)^2} 1_{[|\hat{p}_{ng_n}-p|>\psi_n]}\frac{1}{2\pi}\infint \frac{|\phi_{emp}(t)\phi_w(h_nt)|^2}{|\phi_Z(t)|^2} dt \right]\\
& \lesssim \frac{1}{\epsilon_n^2}\frac{1}{h_n}\infint \frac{|\phi_w(t)|^2}{|\phi_Z(t/h_n)|^2} dt\operatorname{P}(|\hat{p}_{ng_n}-p|>\psi_n).
\end{align*}
Hence
\begin{equation*}
T_{14}\lesssim \frac{1}{\epsilon_n^2}
\frac{1}{h_n^{2\beta+1}}\operatorname{P}(|\hat{p}_{ng_n}-p|>\psi_n)
\end{equation*}
in the ordinary smooth case, and
\begin{equation*}
T_{14}\lesssim
\begin{cases}
\frac{1}{\epsilon_n^2} \frac{1}{h_n}e^{2/(\gamma h_n^{\beta})}\operatorname{P}(|\hat{p}_{ng_n}-p|>\psi_n), & \mbox{if } \beta_0\geq 0,\\
\frac{1}{\epsilon_n^2} h_n^{2\beta_0-1} e^{2/(\gamma h_n^{\beta})}\operatorname{P}(|\hat{p}_{ng_n}-p|>\psi_n), & \mbox{if } \beta_0<0
\end{cases}
\end{equation*}
in the supersmooth case, cf.\ pp.\ 1265-1266 of \citet{fan1}. Similar order bounds are true for $\sup_{p,f} T_{14}.$ Again by Lemma \ref{lemma3} and our conditions on $h_n$ and $\epsilon_n,$ we have $\sup_{p,f} T_{14} \lesssim h_n^{2\alpha}.$

To complete establishing
an upper bound on $T_{10},$ it remains to study $T_{12}.$ As in the case of $T_{11},$ by the
$c_2$-inequality
\begin{align*}
T_{12} & \lesssim \frac{1}{h_n}\infint (w(x))^2dx\ex\left[\frac{(\hat{p}_{ng_n}-p)^2}{(1-\hat{p}_{ng_n})^2(1-p)^2} 1_{[|\hat{p}_{ng_n}-p|\leq\psi_n]} \right]\\
&+\infint\ex\left[(\hat{f}_{nh_n}(x))^2\frac{(\hat{p}_{ng_n}-p)^2}{(1-\hat{p}_{ng_n})^2(1-p)^2} 1_{[|\hat{p}_{ng_n}-p|\leq\psi_n]}\right]dx
\end{align*}
holds. By Lemma \ref{lemma1} of the Appendix the first term on the righthand side is up to a constant bounded by $(1/h_n)g_n^{2\alpha+1}$ and hence is of order $g_n^{2\alpha}.$ The same is true for its supremum over $p$ and $f.$ As far as the second term is concerned, in the supersmooth case it is bounded by
$\psi_n^2\infint\ex[(\hat{f}_{nh_n}(x))^2]dx.$ It follows from the upper bounds on $\sup_{p,f}T_7$ and $\sup_{p,f}T_8$ that in the supersmooth case we have $\sup_{p,f} T_{12}\lesssim h_n^{2\alpha}.$ As far as the ordinary smooth case is concerned,
\begin{multline*}
\infint\ex\left[(\hat{f}_{nh_n}(x))^2\frac{(\hat{p}_{ng_n}-p)^2}{(1-\hat{p}_{ng_n})^2(1-p)^2} 1_{[|\hat{p}_{ng_n}-p|\leq\psi_n]}\right]dx\\
= \infint\ex\left[(\hat{f}_{nh_n}(x))^2\right] dx \, \ex\left[ \frac{(\hat{p}_{ng_n}-p)^2}{(1-\hat{p}_{ng_n})^2(1-p)^2} 1_{[|\hat{p}_{ng_n}-p|\leq\psi_n]}\right]
\end{multline*}
holds. This is precisely the place where we use independence between $\hat{f}_{nh}(x)$ and $\hat{p}_{ng_n}$ implied by sample splitting, cf.\ the remarks around \eqref{ssplf}. Then in this case too $\sup_{p,f} T_{12}\lesssim h_n^{2\alpha}.$ Had not we used the sample splitting trick, in the above display we would have to apply the Cauchy-Schwarz inequality apparently leading to rather lengthy computations.

Combination of the bounds on $\sup_{p,f}T_{11}$ and $\sup_{p,f}T_{12}$ implies that $\sup_{f,p}T_{10}\lesssim
h_n^{2\alpha}.$ The bounds on and $\sup_{f,p}T_{9}$ and $\sup_{f,p}T_{10}$ induce the bound on $T_2.$ The statement of the theorem then
follows from the bounds on $T_1$ and $T_2.$
\end{proof}

\begin{proof}[Proof of Theorem \ref{lbnd-f}] The result is a straightforward consequence of Theorem 2.14 of \cite{meister}.
\end{proof}

\begin{proof}[Proof of Theorem \ref{lbnd-p}]
A general idea of the proof can be outlined as
follows: we will consider two pairs
$(p_1,f_1)$ and $(p_2,f_2)$ (depending
on $n$) of the parameter $(p,f)$ that parametrises the density of $X,$ such that the probabilities $p_1$ and $p_2$
are separated as much as possible, while at the same time the
corresponding product densities $q_1^{\otimes n}$ and
$q_2^{\otimes n}$ of observations $X_1,\ldots,X_n$ are close in
the $\chi^2$-divergence and hence cannot be distinguished well
using the observations $X_1,\ldots,X_n.$ By Lemma 8 of \cite{tsyb1} the squared distance
between $p_1$ and $p_2$ will then
give (up to a constant that does not depend on $n$) the desired lower bound \eqref{lbndp} for estimation of
$p.$

Our construction of the two alternatives $(p_1,f_1)$ and $(p_2,f_2)$ is partially motivated by the construction used in the proof of Theorem 3.5 of \cite{chen}. Let $\lambda_1=\lambda+\delta^{\alpha+1/2},$ where $\lambda>0$ is a fixed constant and $\delta\downarrow 0$ as $n\rightarrow\infty.$ Define $p_1=e^{-\lambda_1}$ and notice that $p_1\in[0,1).$ Next set $\phi_{g_1}(t)=e^{-|t|}$ and observe that this is the characteristic function corresponding to the Cauchy density $g_1(x)=1/(\pi(1+x^2)).$ Finally, define
\begin{equation*}
\phi_{f_1}(t)=\frac{1}{e^{\lambda_1}-1}\left( e^{\lambda_1\phi_{g_1}(t)} - 1 \right).
\end{equation*}
 Denote by $W_j$ the i.i.d.\ random variables that have the common density $g_1$ and by $N_{\lambda_1}$ the random variable that has Poisson distribution with parameter $\lambda_1.$ Then the function $\phi_{f_1}$ will be the characteristic function corresponding to the density $f_1$ of the Poisson sum $Y=\sum_{j=1}^{N_{\lambda_1}}W_j$ of i.i.d.\ $W_j$'s conditional on the fact that the number of its summands $N_{\lambda_1}>0,$ see pp.\ 14--15 of \cite{thesis}. Notice that we have an inequality
\begin{equation*}
|\phi_{f_1}(t)|\leq \frac{\lambda_1 e^{\lambda_1}}{e^{\lambda_1}-1}|\phi_{g_1}(t)|,
\end{equation*}
cf.\ inequality (2.10) on p.\ 22 of \cite{thesis}. Keeping this inequality in mind, without loss of generality we can assume that $K_{\Sigma}$ is already such that $\phi_{f_1}\in\Sigma(\alpha,K_{\Sigma}/4).$ Otherwise we can always consider $\phi_{g_1}(t)=e^{-\alpha^{\prime}|t|}$ with a fixed and large enough constant $\alpha^{\prime}>0,$ so that $\phi_{f_1}\in\Sigma(\alpha,K_{\Sigma}/4).$ It is not difficult to see that the fact that $\alpha^{\prime}\neq 1$ will not affect seriously our subsequent argumentation in this proof. Next define the density $q_1$ corresponding to the pair $(p_1,f_1)$ via its characteristic function
\begin{equation*}
\phi_{q_1}(t)=(p_1+(1-p_1)\phi_{g_1}(t))e^{-t^2/2}
\end{equation*}
and remark that it has the convolution structure required for our problem.

Now we proceed to the definition of the second alternative $(p_2,f_2).$ Set $\lambda_2=\lambda$ and $p_2=e^{-\lambda_2}.$ The fact that $p_2\in[0,1)$ follows from the fact that $\lambda>0.$ Let $H$ be a function, such that its Fourier transform $\phi_H$ is symmetric and real-valued with support on $[-2,2],$ $\phi_H(t)=1$ for $t\in[-1,1]$ and $\phi_H$ is two times continuously differentiable. Such a function can be constructed e.g.\ in the same way as a flat-top kernel in Section 3 of \cite{politis}. Define
\begin{equation*}
\phi_{g_2}(t)=\phi_{g_1}(t)+\tau(t),
\end{equation*}
where the perturbation function $\tau$ is given by
\begin{equation*}
\tau(t)=\frac{\delta^{\alpha+1/2}}{\lambda_2}(\phi_{f_1}(t)-1)\phi_H(\delta t).
\end{equation*}
We claim that for all $n$ large enough $\phi_{g_2}$ is a characteristic function, i.e.\ its inverse Fourier transform $g_2$ is a probability density. This involves showing that $g_2$ integrates to one and is nonnegative. The former easily follows from the fact that
\begin{equation}
\label{intone}
\infint g_2(x)dx=\phi_{g_2}(0)=\phi_{g_1}(0)=1,
\end{equation}
since $\tau(0)=0$ by construction and $\phi_{g_1}$ is a characteristic function. As far as the latter is concerned, we argue as follows: observe that $g_2$ is real-valued, because $\phi_{g_2}$ is symmetric and real-valued. By the Fourier inversion argument
\begin{equation*}
\sup_x |g_2(x)-g_1(x)|\leq \frac{1}{2\pi} \infint|\tau(t)|dt\rightarrow 0
\end{equation*}
as $n\rightarrow\infty,$ by definition of $\tau$ and because $\delta\rightarrow 0.$ Since $g_1,$ being the Cauchy density, is strictly positive on the whole real line, provided $n$ is large enough it follows that
\begin{equation}
\label{g2positive}
g_2(x)\geq 0, \quad x\in B,
\end{equation}
where $B$ is a certain neighbourhood around zero. Next, we need to consider those $x$'s, that lie outside this certain fixed neighbourhood of zero. We have
\begin{align*}
g_2(x)&=\frac{1}{2\pi}\infint e^{-itx} \left( \phi_{g_1}(t) + \frac{\delta^{\alpha+1/2}}{\lambda_2}(\phi_{g_1}(t)-1)\phi_H(\delta t) \right)dt\\
&=\frac{1}{2\pi}\infint e^{-itx} \left( \left(1+\frac{\delta^{\alpha+1/2}}{\lambda_2}\right)\phi_{g_1}(t) - \frac{\delta^{\alpha+1/2}}{\lambda_2}\phi_{g_1}(t) + \frac{\delta^{\alpha+1/2}}{\lambda_2}(\phi_{g_1}(t)-1)\phi_H(\delta t) \right)dt\\
&= \left(1+\frac{\delta^{\alpha+1/2}}{\lambda_2}\right) g_1(x)+\frac{\delta^{\alpha+1/2}}{\lambda_2}\frac{1}{2\pi}\infint e^{-itx}\phi_{g_1}(t)( \phi_H(\delta t) -1 )dt\\
&-\frac{\delta^{\alpha+1/2}}{\lambda_2}\frac{1}{2\pi}\infint e^{-itx} \phi_H(\delta t)dt\\
&=T_1(x)+T_2(x)+T_3(x).
\end{align*}
Both $T_2(x)$ and $T_3(x)$ are real-valued by symmetry of $\phi_{g_1}$ and $\phi_H$ and the fact that these Fourier transforms are real-valued. Consequently, $g_2$ itself is also real-valued. Since $g_1$ is the Cauchy density and $\delta>0,$ the inequality
\begin{equation}
\label{t1x}
T_1(x)\geq \frac{1}{\pi}\frac{1}{1+x^2}
\end{equation}
holds for all $x\in\mathbb{R}.$ Assuming that $x\neq 0$ and integrating by parts, we get
\begin{align*}
T_2(x)&= -\frac{1}{ix} \frac{\delta^{\alpha+1/2}}{\lambda_2}  \frac{1}{2\pi}\int_{\mathbb{R}\setminus[-\delta^{-1},\delta^{-1}]} \phi_{g_1}(t)( \phi_H(\delta t) -1 ) d e^{-itx}\\
& = \frac{1}{ix} \frac{\delta^{\alpha+1/2}}{\lambda_2}  \frac{1}{2\pi} \int_{\mathbb{R}\setminus[-\delta^{-1},\delta^{-1}]} e^{-itx} [ \phi_{g_1}(t)( \phi_H(\delta t) -1 ) ]^{\prime}dt.
\end{align*}
Applying integration by parts to the last equality one more time, we obtain that
\begin{equation*}
T_2(x) =  \frac{1}{x^2} \frac{\delta^{\alpha+1/2}}{\lambda_2}\frac{1}{2\pi} \int_{\mathbb{R}\setminus[-\delta^{-1},\delta^{-1}]} e^{-itx} [ \phi_{g_1}(t)( \phi_H(\delta t) -1 ) ]^{\prime\prime}dt,
\end{equation*}
which implies that
\begin{equation*}
|T_2(x)|\leq  \frac{1}{x^2} C \delta^{\alpha+1/2} \int_{\mathbb{R}\setminus [-\delta^{-1},\delta^{-1}]} | [ \phi_{g_1}(t)( \phi_H(\delta t) -1 ) ]^{\prime\prime} | dt,
\end{equation*}
where the constant $C$ does not depend on $x$ and $n.$
Since $\delta\rightarrow 0$ and the first and the second derivatives of $\phi_H$ are bounded on $\mathbb{R},$ it follows that
\begin{equation*}
|T_2(x)|\leq \frac{1}{x^2} C^{\prime} \delta^{\alpha+1/2}\int_{t>\delta^{-1}} e^{-t}dt,
\end{equation*}
where the constant $C^{\prime}$ is independent of $n$ and $x.$  In particular,
\begin{equation}
\label{t2x}
|T_2(x)|\leq C^{\prime}\delta^{\alpha+1/2} \frac{1}{x^2}
\end{equation}
for all $n$ large enough. Finally, using integration by parts twice, one can also show that for $x\neq 0$
\begin{equation*}
T_3(x)=\frac{1}{x^2}\frac{\delta^{\alpha+5/2}}{\lambda_2}\frac{1}{2\pi}\infint e^{-itx}\phi_H^{\prime\prime}(\delta t)dt
\end{equation*}
and hence
\begin{equation}
\label{t3x}
|T_3(x)|\leq C^{\prime\prime} \delta^{\alpha+3/2}\frac{1}{x^2},
\end{equation}
where the constant $C^{\prime\prime}$ does not depend on $n$ and $x.$ Therefore, by gathering \eqref{t1x}--\eqref{t3x}, we conclude for all $n$ large enough and all $x\in\mathbb{R}$ the inequality
\begin{equation*}
g_2(x)=T_1(x)+T_2(x)+T_3(x) \geq 0
\end{equation*}
is valid. Combining this with \eqref{intone}, we obtain that $g_2$ is a probability density.

Now we turn to the model defined by the pair $(p_2,f_2).$ Again by the argument on pp.\ 22--23 of \cite{thesis},
\begin{equation*}
|\phi_{f_2}(t)|\leq \frac{\lambda_2 e^{\lambda_2}}{e^{\lambda_2}-1}|\phi_{g_2}(t)|.
\end{equation*}
Notice that by selecting $\alpha^{\prime}$ in the definition of $\phi_{g_1}(t)=e^{-\alpha^{\prime}|t|}$ large enough and $\lambda$ large enough, one can arrange that $f_2\in\Sigma(\alpha,K_{\Sigma}),$ at least for all $n$ large enough. Without loss of generality we take $\alpha^{\prime}=1.$ Set
\begin{equation*}
\phi_{q_2}(t)=(p_2+(1-p_2)\phi_{g_2}(t))e^{-t^2/2}.
\end{equation*}
This has the convolution structure as needed in our problem. Hence both pairs $(p_1,f_1)$ and $(p_2,f_2)$ belong to the class required in the statement of the theorem and  generate the required models.

It is easy to see that
\begin{equation}
\label{porder}
|p_2-p_1|\asymp \delta^{\alpha+1/2}
\end{equation}
as $\delta\rightarrow 0,$ where $\asymp$ means that two sequences are asymptotically of the same order. Consequently, by Lemma 8 of \cite{tsyb1} the lower bound in \eqref{lbndp} will be of order $\delta^{2\alpha+1},$ provided we can prove that $n\chi^2(q_2,q_1)\rightarrow 0$
as $n\rightarrow\infty$  for an appropriate $\delta\rightarrow 0.$ Here $\chi^2(q_2,q_1)$ is the $\chi^2$ divergence between the probability measures with densities $q_2$ and $q_1,$ i.e.\
\begin{equation*}
\chi^2(q_2,q_1) = \infint \frac{(q_2(x)-q_1(x))^2}{q_1(x)}dx,
\end{equation*}
see p.\ 86 in \cite{tsyb}.

Notice that we have
\begin{equation*}
q_1(x)=e^{-\lambda_1}\varphi(x)+(1-e^{-\lambda_1})f_1\ast \varphi(x),
\end{equation*}
where $\varphi$ denotes the standard normal density. Let $\delta_1$ denote the first element of the sequence $\delta=\delta_n\downarrow 0.$ Then
\begin{align*}
f_1(x)&=\sum_{n=1}^{\infty} g_1^{\ast n}(x)P(N_{\lambda_1}=n|N_{\lambda_1}>0)\\
& \geq g_1(x) P(N_{\lambda_1}=1|N_{\lambda_1}>0)\\
& = g_1(x) \frac{P(N_{\lambda_1}=1)}{1-P(N_{\lambda_1}=0)}\\
& \geq \frac{\lambda e^{-\lambda-\delta_1^{\alpha+1/2}}}{1-e^{-\lambda_1}} g_1(x),
\end{align*}
cf.\ p.\ 23 in \cite{thesis}. It follows that for all $x$
\begin{equation}
\label{q1ineq}
q_1(x) \geq (1-e^{-\lambda_1}) f_1 \ast \varphi(x) \geq \kappa_A \lambda e^{-\lambda-\delta_1^{\alpha+1/2}} g_1(|x|+A)=c_{\lambda} g_1(|x|+A)
\end{equation}
for some large enough (but fixed) constant $A>0.$ Here the constant $\kappa_A=\int_{-A}^A k(t)dt.$ The inequalities in \eqref{q1ineq} hold, because
\begin{align*}
(1-e^{-\lambda_1})f_1\ast \varphi(x)&=(1-e^{-\lambda_1})\infint f_1(x-t)\varphi(t)dt\\
&\geq \lambda e^{-\lambda-\delta_1^{\alpha+1/2}} \infint g_1(x-t)\varphi(t)dt\\
&\geq \lambda e^{-\lambda-\delta_1^{\alpha+1/2}} \int_{-A}^A g_1(x-t)\varphi(t)dt\\
&\geq g_1(|x|+A) \lambda e^{-\lambda-\delta_1^{\alpha+1/2}} \kappa_A
\end{align*}
by positivity of $g_1$ and $k$ and the fact that the Cauchy density is symmetric at zero and is decreasing on $[0,\infty).$

Now we will use \eqref{q1ineq} to bound the $\chi^2$-divergence between the densities $q_2$ and $q_1.$ Write
\begin{align*}
\chi^2(q_2,q_1)&=\infint\frac{(q_2(x)-q_1(x))^2}{q_1(x)}dx\\
&= \int_{-A}^A\frac{(q_2(x)-q_1(x))^2}{q_1(x)}dx + \int_{\mathbb{R}\setminus[-A,A]}\frac{(q_2(x)-q_1(x))^2}{q_1(x)}dx\\
&=S_1+S_2.
\end{align*}
Using \eqref{q1ineq}, for $S_1$ we have
\begin{equation*}
S_1\leq \frac{1}{c_{\lambda}\inf_{|x|\leq A} g_1(x) } \infint (q_2(x)-q_1(x))^2 dx=c_{\lambda,g_1} \infint (q_2(x)-q_1(x))^2 dx,
\end{equation*}
where $c_{\lambda,g_1}>0$ is a constant.
By Parseval's identity the asymptotic behaviour of the integral on the righthand side of the last equality can be studied as follows,
\begin{align*}
 \infint (q_2(x)-q_1(x))^2 dx & = \frac{1}{2\pi}  \infint |\phi_{q_2}(t)-\phi_{q_1}(t)|^2 dt\\
& = \frac{1}{2\pi} \int_{\mathbb{R}\setminus[-\delta^{-1},\delta^{-1}]} e^{-t^2} \left| e^{\lambda_2(\phi_{g_2}(t)-1)} - e^{\lambda_1(\phi_{g_1}(t)-1)} \right|^2dt\\
& \asymp \frac{1}{2\pi} \int_{\mathbb{R}\setminus[-\delta^{-1},\delta^{-1}]} e^{-t^2}|\delta^{\alpha+1/2}(\phi_{g_1}(t)-1)|^2 |1-\phi_H(\delta t)|^2 dt.
\end{align*}
Using this fact and boundedness of $\phi_H$ on the whole real line, we get that
\begin{equation*}
 \infint (q_2(x)-q_1(x))^2 dx \lesssim \delta^{2\alpha+1} \int_{1/\delta}^{\infty} e^{-t^2}dt\lesssim \delta^{2\alpha+2} e^{-1/\delta^2}.
\end{equation*}
Thus by taking $\delta= c_{\delta}(\log n)^{-1/2}$ with a constant $0<c_{\delta}<1$ we can ensure that the righthand side of the above display is $o(n^{-1})$ and consequently also that $S_1=o(n^{-1}).$

Next we deal with $S_2.$ By \eqref{q1ineq} and Parseval's identity we have that
\begin{equation*}
q_1(x)\geq \frac{c_{\lambda}}{\pi}\frac{1}{1+(|x|+A)^2}.
\end{equation*}
Therefore by Parseval's identity
\begin{equation*}
 S_2\lesssim \int_{\mathbb{R}\setminus[-\delta^{-1},\delta^{-1}]} |[\phi_{q_2}(t)-\phi_{q_1}(t)]^{\prime}|^2dt + \int_{\mathbb{R}\setminus[-\delta^{-1},\delta^{-1}]} |\phi_{q_2}(t)-\phi_{q_1}(t)|^2dt.
\end{equation*}
Exactly by the same type of an argument as for $S_1,$ after some laborious but easy computations, one can show that $S_2=o(n^{-1}),$ provided $\delta\asymp (\log n)^{-1/2}$ with a small enough constant. Consequently, with such a choice of $\delta,$ we have $n\chi^2(q_2,q_1)\rightarrow 0$ as $n\rightarrow \infty$ and the theorem follows from Lemma 8 of \cite{tsyb1} and \eqref{porder}.
\end{proof}

\begin{proof}[Proof of Theorem \ref{lbnd-p-os} ]
We use the same alternatives $(p_1,f_1)$ and $(p_2,f_2)$ as in the proof of Theorem \ref{lbnd-p}. One needs to show that the $\chi^2$-divergence between the corresponding probability densities $q_1$ and $q_2$ is of order $O(n^{-1}).$ The arguments used in the proof of Theorem \ref{lbnd-p} go through and for that end it suffices to show that
\begin{equation}
\label{os1}
\int_{-\infty}^{\infty} | \phi_{q_1}(t) - \phi_{q_2}(t) |^2dt = O(n^{-1})
\end{equation}
and that
\begin{equation}
\label{os2}
\int_{-\infty}^{\infty} | (\phi_{q_1}(t) - \phi_{q_2}(t))^{\prime} |^2dt = O(n^{-1}).
\end{equation}
Observe that for these two integrals to be finite, we need that $\beta>1/2,$ cf.\ the argument below. We have
\begin{align*}
\infint |\phi_{q_2}(t)-\phi_{q_1}(t)|^2 dt
& =  \int_{\mathbb{R}\setminus[-\delta^{-1},\delta^{-1}]} |\phi_Z(t)|^2 \left| e^{\lambda_2(\phi_{g_2}(t)-1)} - e^{\lambda_1(\phi_{g_1}(t)-1)} \right|^2dt\\
& \asymp  \int_{\mathbb{R}\setminus[-\delta^{-1},\delta^{-1}]} |\phi_Z(t)|^2 |\delta^{\alpha+1/2}(\phi_{g_1}(t)-1)|^2 |1-\phi_H(\delta t)|^2 dt.
\end{align*}
Now change the integration variable in the last equality from $t$ to $s=\delta_n t$ and use the fact that for all $s\geq 1$ and for $\delta_n$ small enough by assumption on $\phi_Z$ it holds that
$|\phi_Z(s/\delta_n)| |s/\delta_n|^{\beta}\leq d_1,$
to conclude that the lefthand side of \eqref{os1} is of order $\delta_n^{2\alpha+2\beta}.$ Selecting $\delta_n\asymp n^{-1/(2\alpha+2\beta)}$  then yields \eqref{os1}. A similar argument works in case of \eqref{os2}. We also remark that the condition on $\phi_Z^{\prime}$ given in the statement of the theorem is needed to treat \eqref{os2}. Application of Lemma 8 of \cite{tsyb1} as in Theorem \ref{lbnd-p} concludes the proof.
\end{proof}

\appendix
\section*{Appendix A}

\begin{lem}
\label{lemma1} Let $p^*<1$ and let $\hat{p}_{ng_n}$ be defined by \eqref{ptrunc} (with $p_{ng_n}$ defined by \eqref{png}). Under the same conditions as in Theorem
\ref{thmpng} (i), we have
\begin{equation*}
\sup_{f\in\Sigma(\alpha,K_{\Sigma}),p\in[0,p^*]}\ex[(\hat{p}_{ng_n}-p)^2]\lesssim
n^{-(2\alpha+1)/(2\alpha+2\beta)},
\end{equation*}
while under conditions of Theorem \ref{thmpng} (ii) the inequality
\begin{equation*}
\sup_{f\in\Sigma(\alpha,K_{\Sigma}),p\in[0,p^*]}\ex[(\hat{p}_{ng_n}-p)^2]\lesssim
(\log n)^{-(2\alpha+1)/\beta}
\end{equation*}
holds.
\end{lem}
\begin{proof}[Proof of Lemma \ref{lemma1}]
Introduce the notation
$\sup_{f,p}\equiv\sup_{f\in\Sigma(\alpha,K_{\Sigma}),p\in[0,p^*]}.$ Let $n$ be so large that $p^{\ast}<1-\epsilon_n,$ which is possible, because $p^{\ast}<1$ and $\epsilon_n\downarrow 0.$ Then
\begin{equation*}
\ex[(\hat{p}_{ng_n}-p)^2] \leq \ex[({p}_{ng_n}-p)^2].
\end{equation*}
This and Theorem \ref{thmpng} entail the desired result.
\end{proof}
\begin{lem}
\label{lemma2} Under the same conditions as in Theorem \ref{thmpng} and provided $\epsilon_n=(\log 3n)^{-1},$ the inequality
\begin{equation*}
\sup_{f\in\Sigma(\alpha,K_{\Sigma}),p\in[0,p^*]}\ex\left[ \frac{(\hat{p}_{ng_n}-p)^2}{(1-\hat{p}_{ng_n})^2(1-p)^2}\right] \lesssim g_{n}^{2\alpha+1}
\end{equation*}
holds.
\end{lem}
\begin{proof}
Introduce the sequence
\begin{equation}
\label{psin}
\psi_n=100\sqrt{K_{\Sigma}}U \left\{ \left( \frac{4}{\gamma} \right)^{1/\beta} (\log n)^{-1/\beta} \right\}^{\alpha+1/2}
\end{equation}
and notice that $\psi_n=100\sqrt{K_{\Sigma}}U h_n^{\alpha+1/2}$
in the supersmooth case, i.e in the setting of Theorem \ref{thmpng} (ii). The constants in the definition of $\psi_n$ are rather arbitrary, but they suffice for our purposes. Notice that on the set $\{ |\hat{p}_{ng_n} - p| \leq \psi_n \}$ for all $n$ large enough the inequality
\begin{equation*}
|1-\hat{p}_{ng_n}|\geq 1 - p^{*} - \psi_n
\end{equation*}
holds, because ${\psi}_n\rightarrow 0.$ We have
\begin{align*}
\ex\left[ \frac{(\hat{p}_{ng_n}-p)^2}{(1-\hat{p}_{ng_n})^2(1-p)^2}\right] & = \ex\left[ \frac{(\hat{p}_{ng_n}-p)^2}{(1-\hat{p}_{ng_n})^2(1-p)^2} 1_{[|\hat{p}_{ng_n} - p| \leq \psi_n]}\right]\\
&+\ex\left[ \frac{(\hat{p}_{ng_n}-p)^2}{(1-\hat{p}_{ng_n})^2(1-p)^2} 1_{[|\hat{p}_{ng_n} - p| > \psi_n]}\right]\\
&\lesssim \ex[(\hat{p}_{ng_n}-p)^2]\\
&+\frac{1}{\epsilon_n^2} P( |\hat{p}_{ng_n} - p| > \psi_n )\\
&\lesssim g_n^{2\alpha+1}\\
&+\frac{1}{\epsilon_n^2} P( |\hat{p}_{ng_n} - p| > \psi_n ),
\end{align*}
where in the last inequality we used Lemma \ref{lemma1} and Theorem \ref{thmpng}. It is easy to see that for all $f\in\Sigma(\alpha,K_{\Sigma})$ and $p\in[0,p^*]$ the constants in this chain of inequalities can be made independent of a particular $f$ and a particular $p.$ Then applying Lemma \ref{lemma3} and taking supremum over $f\in\Sigma(\alpha,K_{\Sigma})$ and $p\in[0,p^*]$ on the righthand side of the last equality  establishes the desired result, because
\begin{equation*}
\sup_{p,f}\left(\frac{1}{\epsilon_n^2} P( |\hat{p}_{ng_n} - p| > \psi_n ) \right) = o (g_n^{2\alpha+1})
\end{equation*}
holds under our conditions on $\epsilon_n$ and $g_n.$
\end{proof}
\begin{lem}
\label{lemma3}
Define the sequence $\psi_n$ by \eqref{psin} and let $\epsilon_n=(\log 3n)^{-1}.$ Let $\hat{p}_{ng_n}$ be defined by \eqref{ptrunc} (with $p_{ng_n}$ defined by \eqref{png}). Under the same conditions as in Theorem \ref{thmpng} (i) we have
\begin{multline*}
\sup_{f\in\Sigma(\alpha,K_{\Sigma}),p\in[0,p^*]}P( |\hat{p}_{ng_n} - p| > \psi_n ) \\
\lesssim \frac{1}{\psi_n g_n^{\beta}} \exp\left( - \text{const} \times ng_n^{2\beta}  \right) 
+\exp\left( - \text{const}^{\prime} \times \psi_n^2 g_n^{2\beta} n \right),
\end{multline*}
while under those in Theorem \ref{thmpng} (ii) it holds that
\begin{multline*}
\sup_{f\in\Sigma(\alpha,K_{\Sigma}),p\in[0,p^*]}P( |\hat{p}_{ng_n} - p| > \psi_n ) \\
\lesssim \frac{e^{1/(\gamma g_n^{\beta})}}{\psi_n} \exp\left( - \text{const} \times n e^{-2/(\gamma g_n^{\beta})}  \right)
+\exp\left( - \text{const}^{\prime} \times \psi_n^2 e^{-2/(\gamma g_n^{\beta})} n \right)
\end{multline*}
for the case when $\beta_0\geq 0,$ and
\begin{multline*}
\sup_{f\in\Sigma(\alpha,K_{\Sigma}),p\in[0,p^*]}P( |\hat{p}_{ng_n} - p| > \psi_n ) \\
\lesssim \frac{g_n^{\beta_0} e^{1/(\gamma g_n^{\beta})}}{\psi_n} \exp\left( - \text{const} \times n g_n^{-2\beta_0} e^{-2/(\gamma g_n^{\beta})}  \right) \\
+\exp\left( -\text{const} \times \psi_n^2 g_n^{-2\beta_0} e^{-2/(\gamma g_n^{\beta})} n \right)
\end{multline*}
for the case when $\beta_0<0.$ Here $\text{const}$ and $\text{const}^{\prime}$ are some universal constants (not necessarily the same in all three cases) independent of particular $n,$ $p\in[0,p^*]$ and $f\in\Sigma(\alpha,K_{\Sigma}).$
\end{lem}
\begin{proof}
In this proof we continue numbering of the terms from the proof of Theorem \ref{thmfnhg}, because it is the proof of Theorem \ref{thmfnhg} where this lemma finds its primary use. Observe that
\begin{align*}
\operatorname{P}(|\hat{p}_{ng_n}-p|>\psi_n) & \leq \operatorname{P}(|\ex[\hat{p}_{ng_n}]-p|>\psi_n/2) + \operatorname{P}(|\hat{p}_{ng_n}-\ex[\hat{p}_{ng_n}]|>\psi_n/2)\\
&=T_{15}+T_{16}.
\end{align*}
We have
\begin{align*}
|\ex[\hat{p}_{ng_n}]-p|&\leq |\ex[p_{ng_n}]-p|+ |\ex[\hat{p}_{ng_n}-p_{ng_n}]|\\
&\leq |\ex[p_{ng_n}]-p|+|\ex[(1-\epsilon_n-p_{ng_n})1_{[p_{ng_n}>1-\epsilon_n]}]|\\
&+|\ex[(-1+\epsilon_n-p_{ng_n})1_{[p_{ng_n}<-1+\epsilon_n]}]|\\
&\leq \frac{1}{\sqrt{2}}\sqrt{K_{\Sigma}}Ug_n^{\alpha+1/2}\\
&+\ex[|1-\epsilon_n-p_{ng_n}|1_{[p_{ng_n}>1-\epsilon_n]}]\\
&+\ex[|-1+\epsilon_n-p_{ng_n}|1_{[p_{ng_n}<-1+\epsilon_n]}]\\
&=T_{17}+T_{18}+T_{19}.
\end{align*}
We put the study of $T_{17}$ aside for a while and consider the other two terms. Since $T_{18}$ and $T_{19}$ can be studied in the similar manner, we
consider only $T_{18}.$ Our goal is to show that $T_{18}$ (and by extension $T_{19}$) is negligible in comparison to $T_{17}.$ We have
\begin{equation*}
T_{18}\leq
\left(1+\epsilon_n+\frac{1}{2}\int_{-1}^1\frac{|\phi_u(t)|}{|\phi_Z(t/g_n)|}dt\right)\operatorname{P}(p_{ng_n}>1-\epsilon_n).
\end{equation*}
The righthand side in both cases of the ordinary smooth or
supersmooth $Z$ is of smaller order than $T_{17},$ which can be
seen by employing the arguments on pp.\ 1265-1266 from \citet{fan1} used to bound the integral on the righthand side of the above display
and by the exponential bounds on $\operatorname{P}(p_{ng_n}>1-\epsilon_n),$ which we formulate separately in Lemma \ref{lemma4}. With our conditions on $g_n$ these bounds imply that $\sup_{p,f}T_{18}$ is of lower order than $T_{17}.$ The same is true for $\sup_{p,f}T_{19}.$ As a consequence, $\sup_{p,f} (T_{18}+T_{19})<T_{17}$ for all $n$ large enough. Thus
$T_{15}=0,$ provided $n$ is large enough, because $T_{17}<\psi_n/4$ for all $n$ large enough, and in fact $\sup_{p,f} T_{15}=0$ for all $n$ large enough.

It remains to study $T_{16}.$ We have
\begin{align*}
T_{16} & \leq \operatorname{P}(|\hat{p}_{ng_n}-{p}_{ng_n}|>\psi_n/4)+\operatorname{P}(|{p}_{ng_n}-\ex[\hat{p}_{ng_n}]|>\psi_n/4)\\
& \leq \operatorname{P}(|\hat{p}_{ng_n}-{p}_{ng_n}|>\psi_n/4)+\operatorname{P}(|{p}_{ng_n}-\ex[{p}_{ng_n}]|>\psi_n/8)\\
&+\operatorname{P}(|\ex[{p}_{ng_n}]-\ex[\hat{p}_{ng_n}]|>\psi_n/8)\\
&=T_{20}+T_{21}+T_{22}.
\end{align*}
Notice that
\begin{align*}
T_{20}&\leq \operatorname{P}(|1-\epsilon_n-{p}_{ng_n}| 1_{[p_{ng_n}>1-\epsilon_n]}>\psi_n/8)\\
&+\operatorname{P}(|-1+\epsilon_n-{p}_{ng_n}| 1_{[p_{ng_n}<-1+\epsilon_n]}>\psi_n/8).
\end{align*}
We consider e.g.\ the first term on the righthand side. It is bounded by
\begin{equation*}
\frac{8}{\psi_n}\left( 1-\epsilon_n+\frac{1}{2}\int_{-1}^1 \frac{|\phi_u(t)|}{|\phi_Z(t/g_n)|}dt \right)\operatorname{P}(p_{ng_n}>1-\epsilon_n).
\end{equation*}
Next, as we did above, we use the order bound on the integral on the righthand side, cf.\ pp.\ 1265-1266 in \citet{fan1}, and the exponential bounds on $P(p_{ng_n}>1-\epsilon_n)$ from \eqref{nomeri1} and \eqref{nomeri2} from Lemma \ref{lemma4} to bound the first term in the upper bound on $T_{20}.$ Similar reasoning applies to the second term in the upper bound on $T_{20}.$ There we use Lemma \ref{lemma4half}. These bounds give the first term on the righthand side of the three different formulae in the statement of the lemma.

To bound $T_{21},$ we apply the exponential inequalities from Lemma \ref{lemma5}.
The terms on the righthand side will then give the second terms in the three formulae on the righthand side in the statement of the lemma.

Finally, we turn to $T_{22}.$ Our goal is to
show that there exists $n^{\prime}$ independent of $p$ and $f,$ such that for all $n\geq n^{\prime}$ we have
$T_{22}=0.$ It holds that
\begin{align*}
|\ex[{p}_{ng_n}]-\ex[\hat{p}_{ng_n}]|&\leq \ex[|p_{ng_n}-1+\epsilon_n| 1_{[p_{ng_n}>1-\epsilon_n]}]\\
&+\ex[|p_{ng_n}+1-\epsilon_n| 1_{[p_{ng_n}<1-\epsilon_n]}].
\end{align*}
As the arguments for both terms on the righthand side are similar, we consider only the first term.
We have
\begin{equation*}
\ex[|p_{ng_n}-1+\epsilon_n| 1_{[p_{ng_n}>1-\epsilon_n]}]\leq \left( 1+\epsilon_n+\frac{1}{2}\int_{-1}^1 \frac{|\phi_u(t)|}{|\phi_Z(t/g_n)|}dt \right)\operatorname{P}(p_{ng_n}>1-\epsilon_n).
\end{equation*}
By Lemmas \ref{lemma4} and \ref{lemma4half} and the argument as on pp.\ 1265-1266 of \citet{fan1} the righthand side is negligible compared to $\psi_n$ and it
follows that $T_{22}$ is zero for all large enough $n.$ In fact $n^{\prime}$ can be found, such that this holds true uniformly in $p$ and $f$ for all $n\geq n^{\prime}.$ Gathering all the above bounds entails the statement of the lemma.
\end{proof}

\begin{lem}
\label{lemma4}
Let $p_{ng_n}$ be defined by \eqref{png}. Under the conditions of Theorem \ref{thmpng} (i) we have
\begin{equation}
\label{nomeri1}
\sup_{p\in[0,p^*],f\in\Sigma(\alpha,K_{\Sigma})}\operatorname{P}(p_{ng_n}>1-\epsilon_n)\lesssim \exp\left(
-const\times ng_n^{2\beta}
\right),
\end{equation}
while under conditions of Theorem \ref{thmpng} (ii) we have
\begin{equation}
\label{nomeri2}
\sup_{p\in[0,p^*],f\in\Sigma(\alpha,K_{\Sigma})}\operatorname{P}(p_{ng_n}>1-\epsilon_n)\lesssim
\begin{cases}
\exp\left(
-{const}\times ne^{-2/(\gamma
g_n^{\beta})} \right), & \mbox{if }\beta_0\geq 0,\\
\exp\left(
- {const}\times ng_n^{-2\beta_0}e^{-2/(\gamma
g_n^{\beta})} \right), & \mbox{if }\beta_0< 0.
\end{cases}
\end{equation}
Here $const$ is a universal constant independent of particular $n,p\in[0,p^*]$ and $f\in\Sigma(\alpha,K_{\Sigma}).$
\end{lem}
\begin{proof}
We have
\begin{align*}
\operatorname{P}(p_{ng_n}>1-\epsilon_n)&=\operatorname{P}(p_{ng_n}-\ex[p_{ng_n}]>1-\epsilon_n-\ex[p_{ng_n}])\\
&\leq \operatorname{P}(|p_{ng_n}-\ex[p_{ng_n}]|>1-\epsilon_n-\ex[p_{ng_n}])\\
&=\operatorname{P}\left(\left|\sum_{j=1}^n U_n\left(\frac{-X_j}{g_n}\right) -\ex\left[\sum_{j=1}^n U_n\left(\frac{-X_j}{g_n}\right)\right]\right|>n\frac{(1-\epsilon_n-\ex[p_{ng_n}])}{\pi}\right),
\end{align*}
where
\begin{equation*}
U_n(x)=\frac{1}{2\pi}\int_{-1}^{1}e^{-itx}\frac{\phi_u(t)}{\phi_Z(t/g_n)}dt.
\end{equation*}
Under the conditions of Theorem \ref{thmpng} (i) we have
\begin{equation*}
|U_n(x)|\leq \frac{C}{2\pi} \frac{1}{g_n^{\beta}},
\end{equation*}
while under those of Theorem \ref{thmpng} (ii) the inequality
\begin{equation*}
|U_n(x)|\leq
\begin{cases}
\frac{C^{'}}{2\pi} e^{1/(\gamma
g_n^{\beta})}, & \mbox{if }\beta_0\geq 0,\\
\frac{C^{\prime\prime}}{2\pi} g_n^{\beta_0}e^{1/(\gamma g_n^{\beta})}, &
\mbox{if }\beta_0< 0
\end{cases}
\end{equation*}
holds. Here $C,C^{'}$ and $C^{\prime\prime}$ are some constants independent of $n.$
By \eqref{starp7} we have
\begin{equation}
\label{nomeri9}
|\ex[p_{ng_n}]| \leq |\ex[p_{ng_n}]-p|+p  \leq
p^*+\frac{1}{\sqrt{2}}\sqrt{K_{\Sigma}}Ug_n^{\alpha+1/2}.
\end{equation}
By taking $n_0$ so large that for all $n\geq n_0$
\begin{equation}
\label{nomeri6}
p^*+\frac{1}{\sqrt{2}}\sqrt{K_{\Sigma}}Ug_n^{\alpha+1/2}<1-\epsilon_n
\end{equation}
holds, one can ensure that uniformly in $f$ and $p,$
$1-\epsilon_n-\ex[p_{ng_n}]>0$ for $n\geq n_0.$ Then by Hoeffding's inequality,
see Lemma A.4 on p.\ 198 of \cite{tsyb}, we obtain
\begin{equation*}
\operatorname{P}(p_{ng_n}>1-\epsilon_n)\leq 2\exp\left(
-2\frac{(1-\epsilon_n-\ex[p_{ng_n}])^2}{C^2}ng_n^{2\beta} \right)
\end{equation*}
for the setting of Theorem \ref{thmpng} (i), and
\begin{equation*}
\operatorname{P}(p_{ng_n}>1-\epsilon_n)\leq
\begin{cases}
2\exp\left(
-2\frac{(1-\epsilon_n-\ex[p_{ng_n}])^2}{(C^{'})^2}ne^{-2/(\gamma
g_n^{\beta})} \right), & \mbox{if }\beta_0\geq 0,\\
2\exp\left(
-2\frac{(1-\epsilon_n-\ex[p_{ng_n}])^2}{(C^{\prime\prime})^2}ng_n^{-2\beta_0}e^{-2/(\gamma
g_n^{\beta})} \right), & \mbox{if }\beta_0< 0
\end{cases}
\end{equation*}
for the setting of Theorem \ref{thmpng} (ii).
Since
\begin{equation}
\label{nomeri7}
1-\epsilon_n-\ex[p_{ng_n}]\geq
1-\epsilon_n-p^*-\frac{1}{\sqrt{2}}\sqrt{K_{\Sigma}}Ug_n^{\alpha+1/2}>0
\end{equation}
for all $n$ large enough and uniformly in $f$ and $p,$ see \eqref{nomeri9}, there exists a constant $const$ independent of $n,p\in[0,p^*]$ and $f\in\Sigma(0,K_{\Sigma}),$ such that
\begin{equation*}
\sup_{p,f}\operatorname{P}(p_{ng_n}>1-\epsilon_n)\lesssim \exp\left(
-{const}\times ng_n^{2\beta}
\right)
\end{equation*}
for the setting of Theorem \ref{thmpng} (i), and
\begin{equation*}
\sup_{p,f}\operatorname{P}(p_{ng_n}>1-\epsilon_n)\lesssim
\begin{cases}
\exp\left(
-{const}\times ne^{-2/(\gamma
g_n^{\beta})} \right), & \mbox{if }\beta_0\geq 0,\\
2\exp\left(
- {const}\times ng_n^{-2\beta_0}e^{-2/(\gamma
g_n^{\beta})} \right), & \mbox{if }\beta_0< 0
\end{cases}
\end{equation*}
for the setting of Theorem \ref{thmpng} (ii). This concludes the proof.
\end{proof}

\begin{lem}
\label{lemma4half}
Let $p_{ng_n}$ be defined by \eqref{png}. Under the conditions of Theorem \ref{thmpng} (i) we have
\begin{equation*}
\sup_{p\in[0,p^*],f\in\Sigma(\alpha,K_{\Sigma})}\operatorname{P}(p_{ng_n}<-1+\epsilon_n)\lesssim \exp\left(
-const\times ng_n^{2\beta}
\right),
\end{equation*}
while under conditions of Theorem \ref{thmpng} (ii) we have
\begin{equation*}
\sup_{p\in[0,p^*],f\in\Sigma(\alpha,K_{\Sigma})}\operatorname{P}(p_{ng_n}<-1+\epsilon_n)\lesssim
\begin{cases}
\exp\left(
-{const}\times ne^{-2/(\gamma
g_n^{\beta})} \right), & \mbox{if }\beta_0\geq 0,\\
\exp\left(
- {const}\times ng_n^{-2\beta_0}e^{-2/(\gamma
g_n^{\beta})} \right), & \mbox{if }\beta_0< 0.
\end{cases}
\end{equation*}
Here $const$ is a universal constant independent of particular $n,p\in[0,p^*]$ and $f\in\Sigma(\alpha,K_{\Sigma}).$
\end{lem}

\begin{proof}
The proof is analogous to the proof of Lemma \ref{lemma4} and is therefore omitted.
\end{proof}

\begin{lem}
\label{lemma5}
Let $p_{ng_n}$ be defined by \eqref{png}. Under the conditions of Theorem \ref{thmpng} (i) we have
\begin{equation}
\label{nomeri3}
\sup_{p\in[0,p^*,f\in\Sigma(\alpha,K_{\Sigma})]}P(|p_{ng_n} -\ex[p_{ng_n}]|>\psi_n/8)\lesssim \exp\left( - {const}^{\prime} \times\psi_n^2 ng_n^{2\beta} \right),
\end{equation}
while under conditions of Theorem \ref{thmpng} (ii)
\begin{equation}
\label{nomeri4}
\sup_{p\in[0,p^*,f\in\Sigma(\alpha,K_{\Sigma})]}P(|p_{ng_n} -\ex[p_{ng_n}]|>\psi_n/8)\lesssim \exp\left( - {const}^{\prime} \times\psi_n^2 n e^{2/(\gamma g_n^{\beta})} \right)
\end{equation}
holds. Here $const^{\prime}$ is a universal constant independent of particular $n,p\in[0,p^*]$ and $f\in\Sigma(\alpha,K_{\Sigma}).$
\end{lem}
\begin{proof}
These inequalities can be established by using Hoeffding's inequality in the same way as the exponential bounds on $P(p_{ng_n}>1-\epsilon_n)$ from Lemma \ref{lemma4}.
\end{proof}

\bibliographystyle{plainnat}

\end{document}